\documentclass{amsart}
\numberwithin{equation}{section}
\numberwithin{figure}{section}
\usepackage{amssymb}
\let\cal\mathcal
\def\Ascr{{\cal A}}
\def\Bscr{{\cal B}}
\def\Cscr{{\cal C}}

\def\Escr{{\cal E}}
\def\Fscr{{\cal F}}
\def\Gscr{{\cal G}}

\def\Lscr{{\cal L}}
\def\Mscr{{\cal M}}

\def\Oscr{{\cal O}}

\def\Rscr{{\cal R}}
\def\Sscr{{\cal S}}

\def\Uscr{{\cal U}}
\def\Vscr{{\cal V}}

\let\blb\mathbb
\def\CC{{\blb C}}
\def\DD{{\blb D}}

\def \ZZ{{\blb Z}}

\newcommand{\proj}{\operatorname{proj}}

\def\Id{\operatorname{id}}
\def\pr{\mathop{\text{pr}}\nolimits}

\def\Der{\operatorname{Der}}

\def\Ab{\mathbb{Ab}}

\def\Lotimes{\overset{L}{\otimes}}

\def\Mod{\operatorname{Mod}}

\def\Sk{\operatorname{Sk}}

\def\Ch{\mathop{\mathrm{Ch}}}

\def\Qch{\operatorname{Qch}}
\def\coh{\mathop{\text{\upshape{coh}}}}

\def\Spec{\operatorname {Spec}}

\def\Ext{\operatorname {Ext}}
\def\Hom{\operatorname {Hom}}

\def\RHom{\operatorname {RHom}}

\def\ker{\operatorname {ker}}

\def\Tor{\operatorname {Tor}}

\def\add{\operatorname {add}}

\def\Tot{\operatorname {Tot}}

\def\r{\longrightarrow}

\DeclareMathOperator{\Ind}{Ind}

\DeclareMathOperator{\Aut}{Aut}

\let\dirlim\injlim
\let\invlim\projlim

\newtheorem{lemma}{Lemma}[section]
\newtheorem{proposition}[lemma]{Proposition}

\newtheorem{theorem}[lemma]{Theorem}
\newtheorem{corollary}[lemma]{Corollary}

\newtheorem{lemmas}{Lemma}[subsection]
\newtheorem{propositions}[lemmas]{Proposition}
\newtheorem{propositiondefinitions}[lemmas]{Proposition-Definition}
\newtheorem{theorems}[lemmas]{Theorem}
\newtheorem{corollarys}[lemmas]{Corollary}

\theoremstyle{definition}

\newtheorem{definitions}[lemmas]{Definition}

{

\newtheorem{step}{Step}

}

\theoremstyle{remark}

\newtheorem{remark}[lemma]{Remark}
\newtheorem{remarks}[lemmas]{Remark}

\newdimen\uboxsep \uboxsep=1ex
\def\uboxn#1{\vtop to 0pt{\hrule height 0pt depth 0pt\vskip\uboxsep
\hbox to 0pt{\hss #1\hss}\vss}}

\def\uboxs#1{\vbox to 0pt{\vss\hbox to 0pt{\hss #1\hss}
\vskip\uboxsep\hrule height 0pt depth 0pt}}

\usepackage{amscd}
\usepackage[all]{xy}

\def\redef#1#2{\def#1{#2}}
\redef{\ovl}{\overline}
\redef{\Kern}{\mathrm{ker}}
\redef{\Ob}{\mathrm{Ob}}
\redef{\Mor}{\mathrm{Mor}}
\redef{\Ch}{\mathrm{Ch}}
\redef{\K}{\mathrm{K}}
\redef{\D}{\mathrm{D}}
\redef{\Cokern}{\mathrm{Coker}}
\redef{\Beeld}{\mathrm{Im}}
\redef{\Sub}{\mathsf{Sub}}
\renewcommand{\lim}{\mathrm{lim}}
\redef{\colim}{\mathrm{colim}}
\redef{\lrank}{\mathrm{rank}}
\redef{\cham}{\mathrm{cham}}
\redef{\codim}{\mathrm{codim}}
\redef{\Aut}{\mathrm{Aut}}
\redef{\proj}{\mathrm{proj}}
\redef{\flag}{\mathrm{flag}}
\redef{\type}{\mathrm{type}}
\redef{\Tor}{\mathrm{Tor}}
\redef{\Ext}{\mathrm{Ext}}
\redef{\Hom}{\mathrm{Hom}}
\redef{\Def}{\mathrm{Def}}
\redef{\ddef}{\mathrm{def}}
\redef{\deff}{\mathrm{def}}
\redef{\Sk}{\mathrm{Sk}}
\redef{\pr}{\mathrm{pr}}
\redef{\inn}{\mathrm{in}}
\redef{\es}{\mathrm{es}}
\redef{\op}{^{\mathrm{op}}}
\redef{\Z}{\mathbb{Z}}
\redef{\Q}{\mathbb{Q}}
\redef{\R}{\mathbb{R}}
\redef{\N}{\mathbb{N}}
\redef{\C}{\mathbb{C}}
\redef{\AAA}{\mathfrak{a}}
\redef{\BBB}{\mathfrak{b}}
\redef{\CCC}{\mathfrak{c}}
\redef{\DDD}{\mathfrak{d}}
\redef{\EEE}{\mathfrak{e}}
\redef{\III}{\mathfrak{i}}
\redef{\JJJ}{\mathfrak{j}}
\redef{\KKK}{\mathfrak{k}}
\redef{\OOO}{\mathfrak{o}}
\redef{\PPP}{\mathfrak{p}}
\redef{\QQQ}{\mathfrak{q}}
\redef{\SSS}{\mathfrak{s}}
\redef{\TTT}{\mathfrak{t}}
\redef{\XXX}{\mathfrak{x}}
\redef{\YYY}{\mathfrak{y}}
\redef{\UUU}{\mathfrak{u}}
\redef{\VVV}{\mathfrak{v}}
\redef{\WWW}{\mathfrak{w}}
\redef{\FFF}{\mathfrak{f}}
\redef{\GGGG}{\mathfrak{g}}
\redef{\LLLL}{\mathfrak{l}}
\redef{\Set}{\ensuremath{\mathsf{Set}} }
\redef{\Ab}{\ensuremath{\mathsf{Ab}} }
\redef{\Mod}{\ensuremath{\mathsf{Mod}} }
\redef{\Qch}{\ensuremath{\mathsf{Qch}} }
\redef{\coh}{\ensuremath{\mathsf{coh}} }
\redef{\Dif}{\ensuremath{\mathsf{Dif}} }
\redef{\Pre}{\ensuremath{\mathsf{PreMod}} }
\redef{\Pre}{\ensuremath{\mathsf{Pr}} }
\redef{\Sh}{\ensuremath{\mathsf{Sh}} }
\redef{\mmod}{\ensuremath{\mathsf{mod}} }
\redef{\RMod}{\ensuremath{\mathsf{Mod}(R)} }
\redef{\SMod}{\ensuremath{\mathsf{Mod}(S)} }
\redef{\Rmod}{\ensuremath{\mathsf{mod}(R)} }
\redef{\Smod}{\ensuremath{\mathsf{mod}(S)} }
\redef{\ModR}{\ensuremath{\mathsf{Mod}(R)} }
\redef{\ModS}{\ensuremath{\mathsf{Mod}(S)} }
\redef{\modR}{\ensuremath{\mathsf{mod}(R)} }
\redef{\modS}{\ensuremath{\mathsf{mod}(S)} }

\redef{\AAAMod}{\ensuremath{{\mathsf{Mod}}{(\mathfrak{a})}} }
\redef{\ModAAA}{\ensuremath{{\mathsf{Mod}}{(\mathfrak{a})}} }
\redef{\BBBMod}{\ensuremath{{\mathsf{Mod}}{(\mathfrak{b})}} }
\redef{\ModBBB}{\ensuremath{{\mathsf{Mod}}{(\mathfrak{b})}} }
\redef{\AAAmod}{\ensuremath{{\mathsf{mod}}{(\mathfrak{a})}} }
\redef{\modAAA}{\ensuremath{{\mathsf{mod}}{(\mathfrak{a})}} }
\redef{\BBBmod}{\ensuremath{{\mathsf{mod}}{(\mathfrak{b})}} }
\redef{\modBBB}{\ensuremath{{\mathsf{mod}}{(\mathfrak{b})}} }

\redef{\RCat}{\ensuremath{{\mathsf{Cat}}(R)} }
\redef{\SCat}{\ensuremath{{\mathsf{Cat}}(S)} }
\redef{\CatR}{\ensuremath{{\mathsf{Cat}}(R)} }
\redef{\CatS}{\ensuremath{{\mathsf{Cat}}(S)} }

\redef{\Rng}{\ensuremath{\mathsf{Rng}} }
\redef{\Cat}{\ensuremath{\mathsf{Cat}} }
\redef{\PreCat}{\ensuremath{\mathsf{PreCat}} }
\redef{\Gd}{\ensuremath{\mathsf{Gd}} }
\redef{\PreCX}{\ensuremath{\mathsf{Pre}_{\mathcal{C}}(X)}}
\redef{\PreCY}{\ensuremath{\mathsf{Pre}_{\mathcal{C}}(Y)}}
\redef{\PreoX}{\ensuremath{\mathsf{Pre}(\theta _X)}}
\redef{\PreAbX}{\ensuremath{\mathsf{Pre}_{\Ab}(X)}}
\redef{\PreRngX}{\ensuremath{\mathsf{Pre}_{\Rng}(X)}}
\redef{\SchCX}{\ensuremath{\mathsf{Sch}_{\mathcal{C}}(X)}}
\redef{\SchCY}{\ensuremath{\mathsf{Sch}_{\mathcal{C}}(Y)}}
\redef{\SchAbX}{\ensuremath{\mathsf{Sch}_{\Ab}(X)}}
\redef{\SchRngX}{\ensuremath{\mathsf{Sch}_{\Rng}(X)}}
\redef{\CatX}{\ensuremath{\mathsf{Cat} (X)}}
\redef{\CatY}{\ensuremath{\mathsf{Cat} (Y)}}
\redef{\CatXp}{\ensuremath{\mathsf{Cat} (X,p)}}
\redef{\CatXA}{\ensuremath{\mathsf{Cat} (X,A)}}
\redef{\CatXfU}{\ensuremath{\mathsf{Cat} (X,f(U))}}
\redef{\PreC}{\ensuremath{\mathsf{Pr}(\mathcal{C})}}
\redef{\IndC}{\ensuremath{\mathsf{Ind}(\mathcal{C})}}
\redef{\Ind}{\ensuremath{\mathsf{Ind}}}
\redef{\Pro}{\ensuremath{\mathsf{Pro}}}
\redef{\Fun}{\ensuremath{\mathsf{Fun}}}
\redef{\Add}{\ensuremath{\mathsf{Add}}}
\redef{\add}{\ensuremath{\mathsf{add}}}
\redef{\Inj}{\ensuremath{\mathsf{Inj}}}
\redef{\Proj}{\ensuremath{\mathsf{Proj}}}
\redef{\Cof}{\ensuremath{\mathsf{Cof}}}
\redef{\Open}{\ensuremath{\mathsf{Open}}}
\redef{\inj}{\ensuremath{\mathrm{Inj}}}
\redef{\cof}{\ensuremath{\mathrm{Cof}}}
\redef{\Fp}{\ensuremath{\mathsf{Fp}}}
\redef{\ShCT}{\ensuremath{\mathsf{Sh}(\mathcal{C},\mathcal{T})}}
\redef{\ShCL}{\ensuremath{\mathsf{Sh}(\mathcal{C},\mathcal{L})}}
\redef{\ShCLnul}{\ensuremath{\mathsf{Sh}(\mathcal{C},\mathcal{L}
_0)}}
\redef{\ShCLplus}{\ensuremath{\mathsf{Sh}(\mathcal{C},\mathcal{L}
^{+})}}
\redef{\SepCT}{\ensuremath{\mathsf{Sep}(\mathcal{C},\mathcal{T})}}
\redef{\SepCL}{\ensuremath{\mathsf{Sep}(\mathcal{C},\mathcal{L})}}
\redef{\SepCLnul}{\ensuremath{\mathsf{Sep}(\mathcal{C},\mathcal{L}
_0)}}
\redef{\SepCLplus}{\ensuremath{\mathsf{Sep}(\mathcal{C},\mathcal{L}
^{+})}}
\redef{\lra}{\longrightarrow}
\redef{\ra}{\rightarrow}
\redef{\Stensor}{\ensuremath{(S \otimes_R -)}}
\redef{\HomS}{\ensuremath{\mathrm{Hom}_R(S,-)}}

\newcommand{\aaa}{\ensuremath{\mathcal{A}}}
\newcommand{\bbb}{\ensuremath{\mathcal{B}}}
\newcommand{\ccc}{\ensuremath{\mathcal{C}}}
\newcommand{\ddd}{\ensuremath{\mathcal{D}}}
\newcommand{\eee}{\ensuremath{\mathcal{E}}}

\newcommand{\ooo}{\ensuremath{\mathcal{O}}}

\newcommand{\qqq}{\ensuremath{\mathcal{Q}}}
\newcommand{\rrr}{\ensuremath{\mathcal{R}}}
\newcommand{\sss}{\ensuremath{\mathcal{S}}}

\newcommand{\uuu}{\ensuremath{\mathcal{U}}}

\newcommand{\xxx}{\ensuremath{\mathcal{X}}}

\redef{\CC}{\ensuremath{\mathbf{C}}}
\redef{\DD}{\ensuremath{\mathbf{D}}}

\newcommand{\DGFun}{\ensuremath{\mathsf{DGFun}} }

\newcommand{\Ho}{\ensuremath{\mathrm{Ho}} }
\dedicatory{Dedicated to Michael Artin on the occasion of his 
seventieth birthday.}
\title{Hochschild cohomology of abelian categories and ringed spaces}
\author{Wendy T. Lowen}
\address{Departement DWIS\\ Vrije Universiteit Brussel\\ Pleinlaan
2\\1050 Brussel\\ Belgium}
\email[Wendy T. Lowen]{wlowen@vub.ac.be}
\author{Michel Van den Bergh}
\address{Departement WNI\\Limburgs Universitair
Centrum\\ Universitaire Campus\\ Building D\\ 3590
Diepenbeek\\ Belgium}
\email[Michel Van den Bergh]{vdbergh@luc.ac.be}
\thanks{The first author is an aspirant at the FWO}
\thanks{The second
author is a senior researcher at the FWO}
\keywords{Hochschild cohomology, abelian categories, DG-categories} 
\subjclass{Primary 13D10, 14A22, 18E15} 
\thanks{This paper was completed while both authors were visiting
the Mittag-Leffler institute during the program on non-commutative
geometry in the academic year 2003/2004. We hereby thank the Mittag-Leffler institute for its kind hospitality and for the stimulating
atmosphere it provides.}
\begin{document}
\begin{abstract}
   This paper continues the development of the deformation theory of
   abelian categories introduced in a previous paper by the authors. We
   show first that the deformation theory of abelian categories is
   controlled by an obstruction theory in terms of a suitable notion of
   Hochschild cohomology for abelian categories. We then show that this
   Hochschild cohomology coincides with the one defined by
   Gerstenhaber, Schack and Swan in the case of module categories over
   diagrams and schemes and also with the Hochschild cohomology for
   exact categories introduced recently by Keller. In addition we show  in
   complete generality that Hochschild cohomology satisfies a
   Mayer-Vietoris property and that for constantly ringed spaces it
   coincides with the cohomology of the structure sheaf.
\end{abstract}
\maketitle
\tableofcontents

\section{Introduction}
In the rest of this paper $k$ is an arbitrary commutative base ring
but for simplicity we will assume in this introduction that $k$ is a
field.

\medskip

Motivated by our work on the infinitesimal deformation theory of
abelian categories \cite{lowenvdb1} our aim in this paper is to
develop a theory of Hochschild cohomology for abelian categories and
ringed spaces. The corresponding theory for Hochschild (and cyclic)
\emph{homology} is by now well established
\cite{keller5,loday,weibel}.  The theory for Hochschild cohomology has
a rather different flavour since it is less functorial, but
nevertheless it still has good stability and agreement properties.

\medskip

We start with the case of $k$-linear categories. The Hochschild
complex $\CC(\AAA)$ of a $k$-linear category $\AAA$ is defined by
\begin{equation}
\label{ref-1.1-0}
\CC^p(\AAA)=\prod_{A_0,\dots,A_p\in 
\Ob(\AAA)}\Hom_k(\AAA(A_{p-1},A_p) \otimes_k \dots \otimes_k
\AAA(A_0,A_1), \AAA(A_0,A_p))
\end{equation}
with the usual differential (see \cite{Mi}). As in the case
of associative algebras, the Hochschild complex of $\AAA$ carries a
considerable amount of ``higher structure'' containing in particular
the classical cup-product and the Gerstenhaber bracket. This extra
structure may be summarized conveniently by saying that $\CC(\AAA)$ is
an algebra over the so-called $B_\infty$-operad \cite{getzlerjones,keller6}.

  Let
$\Ascr$ be a $k$-linear  abelian category. In this paper we define the
\emph{Hochschild complex}  of $\Ascr$ as
\begin{equation}
\label{ref-1.2-1}
  \CC_{\mathrm{ab}}(\Ascr) =\CC(\Inj \Ind(\Ascr))
\end{equation}
where $\Ind \Ascr$ is the abelian category of $\Ind$-objects over
$\Ascr$ \cite[Expose I]{SGA4tome1} and $\Inj \Ind(\Ascr)$ denotes the full
subcategory of injective objects in $\Ind(\Ascr)$.  It is understood
here that $\Ind\Ascr$ is computed with respect to a universe in which
$\Ascr$ is small. In the rest of this introduction,  for the purpose 
of exposition,
we will ignore such
settheoretic complications (see
\S\ref{ref-2.1-5},\S\ref{ref-2.5-17} below).

As indicated above the initial motivation behind \eqref{ref-1.2-1} is the
infinitesimal deformation theory of $k$-linear abelian categories
developed in \cite{lowenvdb1} (see \S\ref{ref-3-18} below).  We prove
(Theorem \ref{ref-3.1-20}):
\begin{itemize}
\item the deformation theory of a $k$-linear abelian category $\Ascr$
   is controlled by an obstruction theory involving
   $H\CC^2_{\mathrm{ab}}(\Ascr)$ and $H\CC^3_{\mathrm{ab}}(\Ascr)$.
\end{itemize}
In \S\ref{ref-6-61} we prove basic results about the Hochschild
cohomology of $k$-linear abelian categories. In particular we show
(Theorem \ref{ref-6.6-67} and Corollaries \ref{ref-6.8-69},\ref{ref-6.9-70}):
\begin{itemize}
\item If $\Ascr$ has itself enough injectives then
   $\CC_{\mathrm{ab}}(\Ascr)\cong \CC(\Inj \Ascr))$ (where here and
   below $\cong$ means the existence of an isomorphism in the homotopy
   category of $B_\infty$-algebras). A dual statement holds of course
   if $\Ascr$ has enough projectives.
\item In general we have $\CC_{\mathrm{ab}}(\Ascr)\cong
   \CC_{\mathrm{ab}}(\Ind \Ascr)$.
\item If $A$ is a $k$-algebra then the Hochschild cohomology of the
   abelian category $\Mod(A)$ coincides with the Hochschild cohomology
   of $A$.
\end{itemize}
We also show that the Hochschild complex of an abelian category
$\Ascr$ is the same as the Hochschild complexes of suitable
DG-categories \cite{Keller1} associated to $\Ascr$. The definition of
the Hochschild complex of a DG-category is an obvious extension of
\eqref{ref-1.1-0} (see \S\ref{ref-2.4-13} below).  Let ${}^e\!\!\Ascr$ be
the full DG-subcategory of $C(\Ind(\aaa))$ spanned by all positively
graded complexes of injectives whose only cohomology is in degree zero
and lies in $\aaa$ and let ${}^eD^b(\Ascr)$ be spanned by all left
bounded complexes of injectives with bounded cohomology in~$\aaa$.
Then ${}^e \! D^b(\Ascr)$ is an exact DG-category such that $
H^\ast({}^e \! D^b(\Ascr)) $ is the graded category associated to
$D^b(\Ascr)$. Hence ${}^e \! D^b(\Ascr)$ is a DG-``enhancement''
\cite{BK1,BLL} for the triangulated category $\!\!D^b(\Ascr)$.

We prove (Theorem \ref{ref-6.1-62})
\begin{equation}
\label{ref-1.3-2}
\CC_{\mathrm{ab}}(\Ascr)\cong \CC({}^e\!\!\Ascr)\cong \CC({}^e\!D^b(\Ascr))
\end{equation}
Using \eqref{ref-1.3-2} we may construct a homomorphism of graded rings
(see Prop.\ \ref{ref-4.5-41})
\[
\sigma_\Ascr:H\CC_{\mathrm{ab}}^\ast(\Ascr)\r Z(D^b(\Ascr))
\]
where  we view $D^b(\Ascr)$ as a graded category. The
homogeneous elements of $Z(D^b(\Ascr))$ consist of tuples
$(\phi_M)_{M}$ with $M\in \Ob(D^b(\Ascr))$ and $\phi_M\in
\Ext^\ast_\Ascr(M,M)$ satisfying a suitable compatibility condition
(see \S\ref{ref-4.5-41}). Thus we may think of $\sigma_\Ascr$ as defining
   ``universal'' elements in $\Ext^\ast_\Ascr(M,M)$ for every $M\in
   D^b(\Ascr)$.  These universal elements are closely related to
   Atiyah classes in algebraic geometry. See for
   example~\cite{BuchweitzFlenner}.

   The isomorphisms in \eqref{ref-1.3-2} also connect 
$\CC_{\mathrm{ab}}(\Ascr)$ to
   Keller's recent definition of the Hochschild complex of an exact
   category \cite{keller6}. If $\Escr$ is an exact category then by
   definition $\CC_{\mathrm{ex}}(\Escr)= \CC(\qqq)$ for a DG-quotient
   \cite{keller6,Drinfeld2} $\qqq$ of $Ac^b(\Escr) \lra C^b(\Escr)$, where
   $C^b(\Escr)$ is the DG-category of bounded complexes of
   $\eee$-objects and $Ac^b(\eee)$ is its full DG-subcategory of
   acyclic complexes. If we equip an abelian category $\Ascr$ with its
   canonical exact structure then ${}^eD^b(\Ascr)$ is a DG-quotient of
   $Ac^b(\Escr) \lra C^b(\Escr)$ (see Lemma \ref{ref-6.3-64}).  Hence
   $\CC_{\mathrm{ab}}(\Ascr)\cong \CC_{\mathrm{ex}}(\Ascr)$.

In \S\ref{ref-7-72} we specialize to ringed spaces. If $(X,\Oscr)$
is a ringed space then we define
\[
\CC(X)=\CC_{\mathrm{ab}}(\Mod(X))
\]
where $\Mod(X)$ is the category of sheaves of right $\Oscr$-modules.
Note that in this definition the bimodule structure of $\Oscr$ does
not enter explicitly. We show that $H\CC^\ast(-)$ defines a ``nice''
cohomology theory for $(X,\Oscr)$ in the sense that it has the
following properties (\S\ref{ref-7.4-85} and Theorem \ref{ref-7.9.1-104}).
\begin{itemize}
\item $H\CC^\ast(-)$ is a contravariant functor on open embeddings.
\item Associated to an open covering $X=U\cup V$ there is a
   Mayer-Vietoris sequence
\[
\cdots \rightarrow H\CC^{i-1}(U\cap V)\rightarrow
  H\CC^i(X)\rightarrow H\CC^i(U)\oplus H\CC^i(V)
\rightarrow H\CC^i(U\cap V)\rightarrow \cdots
\]
\item If $\Oscr$ is the constant sheaf $\underline{k}$ with values in $k$
   then
\begin{equation}
\label{ref-1.4-3}
H\CC^\ast(X,\underline{k})\cong H^\ast(X,\underline{k})
\end{equation}
where the righthandside is the usual derived functor cohomology of
$\underline{k}$.
\end{itemize}
  In \cite{Baues} Baues
shows that the singular cochain complex of a topological space is a
$B_\infty$-algebra. Thus \eqref{ref-1.4-3} suggests that
$\CC(X,\underline{\ZZ})$ should be viewed as an algebraic analog of
the singular cochain complex of $X$.

\medskip

We also show that under suitable conditions $H\CC^\ast(X)$ coincides with the
Hochschild cohomology theories for ringed spaces and schemes defined
by Gerstenhaber, Schack and Swan. More precisely we show:
\begin{itemize}
\item Assume that $X$ has a basis $\Bscr$ of \emph{acyclic opens}, i.e.
   for $U \in \bbb$: $H^i(U,\ooo_U) = 0$ for $i>0$.  Then
   $H\CC^\ast(X)$ coincides with the Gerstenhaber-Schack cohomology
   \cite{GS1,GS} of the restriction of $\Oscr$ to $\Bscr$, considered as a
   diagram over the partially ordered set $\Bscr$.
\item If $X$ is a quasi-compact separated scheme then $H\CC^\ast(X)$
   coincides with the Hochschild cohomology for $X$ as defined by Swan
   in \cite{swan}.
\end{itemize}
We recall that for a smooth scheme the Hochschild complex defined
by Swan is quasi-isomorphic to the one defined by Kontsevich \cite{Ko3} in
terms of differential operators (see \cite{swan,Yek2}).

Now let $X$ be a quasi-compact separated scheme and denote by
$\Qch(X)$ the category of quasi-coherent sheaves on $X$. If $X$ is in
addition noetherian then let $\coh(X)$ be the coherent sheaves on $X$.
We prove that there are isomorphisms (see Theorem \ref{ref-7.5.1-88} and
Corollaries \ref{ref-7.7.2-99},\ref{ref-7.7.3-100})
\begin{equation}
\label{ref-1.5-4}
\CC(X)\cong \CC_{\mathrm{ab}}(\Qch(X))\cong \CC_{\mathrm{ab}}(\coh(X))
\end{equation}
whenever the notations make sense.

The first isomorphism in \eqref{ref-1.5-4} is proved by relating the
Hochschild cohomology of $X$ to that of a finite open affine covering
of $X$.  To be more precise let $X=A_1\cup\cdots \cup A_n$ be such a
covering and let $\Ascr$ be the closure of $\{A_1,\ldots,A_n\}$ under
intersections. Since $X$ is separated, $\Ascr$ consists of affine
opens.  We define a linear category $\AAA$ with the same objects as
$\Ascr$ by putting
\[
\AAA(U,V)=\begin{cases}
\Gamma(U,\Oscr_U)&\text{if $U\subset V$}\\
0&\text{otherwise}
\end{cases}
\]
We construct an isomorphism (Corollary
\ref{ref-7.7.2-99})
\[
\CC(X)\cong \CC(\AAA)
\]
In particular if $X=\Spec R$ is itself affine then
$
\CC(X)\cong \CC(R)
$.
\medskip

We are extremely grateful to Bernhard Keller for freely sharing with
us many of his ideas and for making available the preprint
\cite{keller6}. While preparing the current manuscript the authors had
independently discovered the main result of \cite{keller6} (with the same
proof) but nevertheless the presentation in \cite{keller6} made it
possible to clarify and generalize many of our original arguments.

\medskip

The second author learned about the connection
between Hochschild cohomology and Atiyah classes
in an interesting discussion with Ragnar-Buchweitz
  at a sushi restaurant in Berkeley
during the workshop on non-commutative algebraic geometry at MSRI in
Februari 2000.

\section{Preliminaries, conventions and notations}
\subsection{Universes}
\label{ref-2.1-5}
The results in this paper are most conveniently stated for small
categories. However we sometimes need non-small categories as well,
for example to pass from a category to its category of $\Ind$-objects.
Therefore we take the theory of universes as our set theoretical
foundation since this basically allows us to assume that any category
is small. For a brief introduction to the theory of universes and to
some related terminology which we will use in this paper, we refer the
reader to \cite{lowenvdb1}. Our convention is that we fix a universe
$\uuu$, and all terminology (small, Grothendieck, $\dots$) and all
constructions ($\Ab$, $\Mod$, $\Ind$, $\dots$) are implicitly prefixed
by $\uuu$. Unless otherwise specified all categories will be
$\uuu$-categories, i.e. their $\Hom$-sets are $\Uscr$-small.
Individual objects such as rings and modules are also assumed to be
$\Uscr$-small.

\subsection{DG-categories} We will assume that the reader has some
familiarity with DG-categories and model categories. See for example
\cite{Keller1,Hovey}.

Throughout we fix a commutative ring $k$ and we assume that all
categories are $k$-linear. Unadorned tensor products and $\Hom$'s are
over $k$. On first reading one may wish to assume that $k$ is a field
as it technically simplifies many definitions and proofs (see for
example the next section).

Let $\AAA$ be a DG-category. Associated to $\AAA$ is the
corresponding graded category (for which we use no separate notation),
which is obtained by forgetting the differential
and the categories
$H^0(\AAA)$ and $H^\ast(\AAA)$ with
$\Ob(H^0(\AAA))=\Ob(H^\ast(\AAA))=\Ob(\AAA)$ and
\begin{align*}
H^0(\AAA)(A,B)&=H^0(\AAA(A,B))\\
H^\ast(\AAA)(A,B)&=H^\ast(\AAA(A,B))
\end{align*}
$H^0(\AAA)$ is sometimes referred to as the \emph{homotopy category} of $\AAA$.

Now assume that $\AAA$ is small. We consider the DG-category
$$\Dif(\AAA) = \DGFun(\AAA^{\op},C(k))$$ of (right) DG-modules. The
derived category $D(\AAA)$ is the localized category
\[
\Dif(\AAA)[\Sigma^{-1}]
\]
where $\Sigma$ is the class of (pointwise) quasi-isomorphisms
\cite{Keller1}.  In order to work conveniently with $D(\AAA)$ it
is useful to introduce model structures \cite{Hovey} on $\Dif(\AAA)$.

It turns out that $\Dif(\AAA)$ is equipped with two canonical model
structures for which the weak equivalences are the quasi-isomorphisms.
\cite{Hinich}.  For the \emph{projective} model structure the
fibrations are the pointwise epimorphisms and for the
\emph{injective} model structure the cofibrations are the pointwise
monomorphims.

We say $M\in \Dif(\AAA)$ is \emph{fibrant} if $M\r 0$ if a fibration
for the injective model structure and we call $M$ \emph{cofibrant} if
$0\r M$ is a cofibration for the projective model structure.

As usual $D(\AAA)$ is the homotopy category of cofibrant
complexes and also the homotopy category of fibrant complexes and this
makes it easy to construct left and right derived functors.

Let $\frak{b}$ be another small DG-category and let $f: \AAA \lra
\BBB$ be a DG-functor.
$f$ is a
\emph{quasi-equivalence} if
$H^{\ast}(f)$ is fully faithful and $H^0(f)$ is essentially surjective and
it is called a \emph{quasi-isomorphism} if $H^\ast(f)$ is an
isomorphism.

The functor $f$ induces the usual triple of adjoint functors
$(f^\ast,f_\ast,f^!)$ between $\Dif(\AAA)$ and $\Dif(\BBB)$. We denote
the corresponding adjoint functors between $D(\AAA)$ and $D(\BBB)$ by
$(Lf^\ast,f_\ast,Rf^!)$.
\begin{propositions}\cite{Keller1} \label{ref-2.2.1-6}
The functors $(Lf^\ast,f_\ast,Rf^!)$ are equivalences  when $f$ is
a quasi-equivalence.
\end{propositions}

For $M, N \in \Dif(\AAA)$, $\Hom_{\AAA}(M,N)$ denotes $\Dif(\AAA)(M,N)
\in C(k)$.  There is a corresponding derived functor
\[
\RHom_{\AAA}: D(\AAA)^{\op} \times D(\AAA) \lra D(k).
\]
\begin{lemmas} \label{ref-2.2.2-7} If $f$ is a quasi-equivalence then 
the induced map
\[
f_\ast:\RHom_{\AAA}(M,N)\r \RHom_{\BBB}(f_\ast M,f_\ast N)
\]
is a quasi-isomorphism.
\end{lemmas}
\begin{proof}
After replacing $M$ by a cofibrant resolution, $f_\ast$ is defined as
the composition
\[
\RHom_{\AAA}(M,N)=\Hom_{\AAA}(M,N)\r \Hom_{\BBB}(f_\ast M,f_\ast N)
\mathbin{\mathop{\longrightarrow}\limits^{\mathrm{can}}}
\RHom_{\BBB}(f_\ast M,f_\ast N)
\]
and looking at homology we see, using Proposition \ref{ref-2.2.1-6}, that this
is a quasi-isomorphism.
\end{proof}

\subsection{Bimodules and resolutions}
\label{ref-2.3-8}
An $\AAA{-}\BBB$-(DG-)bimodule $X$ is an object of $\Dif(\AAA^{\op} \otimes
  \BBB)$
which  will be denoted by $(B,A) \longmapsto X(B,A)$, where $X(B,A)$ is
contravariant in $B$ and covariant in $A$.
The $\AAA{-}\AAA$-bimodule $(A,A')
\longmapsto \AAA(A,A')$ will be denoted by $\AAA$. For $\AAA{-}\BBB$-bimodules
$X$ and
$Y$, we define the
$\AAA{-}\AAA$ and $\BBB{-}\BBB$-bimodules
\begin{equation}
\label{ref-2.1-9}
\begin{aligned}
\Hom_{\BBB}(X,Y)(A,A') &= \Hom_{\BBB}(X(-,A),Y(-,A'))\\
\Hom_{\AAA^{\op}}(X,Y)(B,B') &= \Hom_{\AAA^{\op}}(X(B',-),Y(B,-)).
\end{aligned}
\end{equation}
For $G \in \Dif(\BBB)$ and $F \in \Dif(\BBB^{\op})$, there is also
a tensor product $G \otimes_{\BBB} F \in C(k)$ with in particular $\BBB(-,B)
\otimes_{\BBB} F = F(B)$ and $G \otimes_{\BBB} \BBB(B,-) = G(B)$ (see for
example
\cite[14.3, 14.4]{Drinfeld2}).  For
$X
\in
\Dif(\AAA^{\op}
\otimes
\BBB)$ and
$Y
\in
\Dif(\BBB^{\op}
\otimes
\CCC)$, we define $X \otimes_{\BBB} Y \in \Dif(\AAA^{\op} \otimes \CCC)$ by
$$(X \otimes_{\BBB} Y)(C,A) = X(-,A) \otimes_{\BBB} Y(C,-)$$
(see also \cite[14.5]{Drinfeld2}).

If $k$ is not a field then in the derived setting the tensor product
$\AAA^{\op} \otimes \BBB$ is philosophically wrong! Instead one should use
something like ``$\AAA^{\op}\Lotimes \BBB$'' but of course this has no
immediate meaning.  Thus when working with bimodules we should assume
that our categories satisfy an appropriate flatness assumption. It
turns out that it is most convenient to assume that our categories are
$k$-cofibrant in the following sense \cite{keller6}:
\begin{definitions}
\label{ref-2.3.1-10}
A DG-category is \emph{$k$-cofibrant} if all its $\Hom$-sets are 
cofibrant in $C(k)$.
\end{definitions}
Recall that if $M\in C(k)$ is cofibrant then it has projective terms
and the functors $\Hom(M,-)$ and $M\otimes-$ preserve acyclic
complexes. This is technically very useful.

\medskip

That Definition \ref{ref-2.3.1-10}  is a good definition follows from
the following result.
\begin{propositiondefinitions}
\label{ref-2.3.2-11}
\begin{enumerate}
\item Let $\AAA$ be a small DG-category. There exists a
   quasi-isomorphism $\ovl{\AAA}\r \AAA$ with $\bar{\AAA}$
   $k$-cofibrant which is surjective on $\Hom$-sets (in the graded
   category). We call such $\ovl{\AAA}\r \AAA$ a \emph{$k$-cofibrant
     resolution} of $\AAA$.
\item If $f:\AAA\r \BBB$ is a DG-functor between small DG-categories
   and $\ovl{\BBB}\r \BBB$ is a $k$-cofibrant resolution of $\BBB$ then
there exists a $k$-cofibrant resolution $\ovl{\AAA}\r \AAA$ together
with a commutative diagram
\[
\begin{CD}
\ovl{\AAA} @>\ovl{f}>> \ovl{\BBB}\\
@VVV @VVV\\
\AAA @>>f> \BBB
\end{CD}
\]
\end{enumerate}
\end{propositiondefinitions}

\begin{proof}
\begin{enumerate}
\item We may take $\bar{\AAA}$ to be a ``semi-free'' resolution of
   $\AAA$. See \cite[Lemma 13.5]{Drinfeld2} and \cite{Hinich}.
\item Again we let $\bar{\AAA}\r\AAA$  be a semi-free resolution. The
   result now follows from \cite[Lemma 13.6]{Drinfeld2}.\qed
\end{enumerate}
\def\qed{}\end{proof} If $\AAA$, $\BBB$ are small DG-categories then
the derived category of $\AAA$-$\BBB$-modules should be
defined as
\begin{equation}
\label{ref-2.2-12}
D(\ovl{\AAA}^{\op} \otimes \ovl{\BBB})
\end{equation}
for $k$-cofibrant resolutions $\ovl{\AAA}\r \AAA$, $\ovl{\BBB}\r
\BBB$. Propositions \ref{ref-2.2.1-6} and \ref{ref-2.3.2-11} insure that this
definition is independent of the chosen resolutions (up to
equivalence). We use this definition in principle, but to make things
not overly abstract we will always indicate the resolutions
$\ovl{\AAA}$ and $\ovl{\BBB}$ explicitly in the notations.
\begin{propositions} Assume that $\AAA$ and $\BBB$ are $k$-cofibrant.
The derived functors
\begin{align*}
   \RHom_{\BBB}:& D(\AAA^{\op} \otimes \BBB)^{\op}
  \times D(\AAA^{\op} \otimes \BBB) \lra D(\AAA^{\op} \otimes \AAA)\\
\RHom_{\AAA^{\op}}:& D(\AAA^{\op} \otimes \BBB)^{\op}
  \times D(\AAA^{\op} \otimes \BBB) \lra D(\BBB^{\op} \otimes \BBB).
\end{align*}
of \eqref{ref-2.1-9} may be computed pointwise  in the sense
\begin{equation}
\begin{aligned}
\RHom_{\BBB}(X,Y)(A,A') &= \RHom_{\BBB}(X(-,A),Y(-,A'))\\
\RHom_{\AAA^\circ}(X,Y)(B,B') &= \RHom_{\AAA^\circ}(X(B',-),Y(B,-)).
\end{aligned}
\end{equation}
\end{propositions}
\begin{proof} Easy.
\end{proof}
\subsection{Hochschild cohomology of DG-categories}
\label{ref-2.4-13}
Let $\AAA$ be a small $k$-cofibrant DG-category and let $M$ be a
$\AAA$-bimodule.
The Hochschild complex
$\CC(\AAA,M)$ of
$\AAA$ with coefficients in an $\AAA{-}\AAA$-bimodule
$M$ is the total complex of the double complex $\DD(\AAA,M)$ with $p$-th column
given by
\begin{equation}
\label{ref-2.4-14}
\prod_{A_0,\dots,A_p}\Hom_k(\AAA(A_{p-1},A_p) \otimes_k \dots \otimes_k
\AAA(A_0,A_1), M(A_0,A_p))
\end{equation}
  and the usual horizontal Hochschild differential.
The Hochschild complex of  $\AAA$ is $\CC(\AAA) = \CC(\AAA,\AAA)$.
There is an isomorphism in $D(k)$
\begin{equation}
\label{ref-2.5-15}
\CC(\AAA,M) \cong \RHom_{\AAA^{\op}
\otimes
\AAA}(\AAA,M).
\end{equation}
The Hochschild complex satisfies a ``limited functoriality''
property. If $j:\AAA\r \BBB$ is a fully faithful map between small
$k$-cofibrant DG-categories then there is an associated map between
Hochschild complexes
\[
j^\ast:\CC(\BBB)\r \CC(\AAA)
\]
given by restricting cocycles. We will usually refer to $j^\ast$ as
the \emph{restriction map}.

It is well-known that $\CC(\AAA)$ carries a considerable amount
of ``higher structure'' containing in particular the classical cup-product and
the Gerstenhaber bracket. This extra structure is important for
deformation theory.  A convenient way of summarizing the extra
structure is by saying that $\CC(\AAA)$ is an algebra over the
$B_\infty$-operad \cite{getzlerjones,keller6} which is an enlargement of the
$A_\infty$-operad.  The map $j^\ast$ introduced above is trivially
compatible with the $B_\infty$-structure.

As the Hochschild complex  involves bimodules,
according to the principles outlined in \S\ref{ref-2.3-8}, its definition
should be modified for non $k$-cofibrant DG-categories. The
appropriate modification was introduced by Shukla and Quillen
\cite{Quillen1,Shukla} in the case of DG-algebras.

Let $\AAA$ be a small DG-category which is not necessarily
$k$-cofibrant and let $M$ be in $\Dif(\AAA^\circ\otimes \AAA)$. We
fix a k-cofibrant resolution $\bar{\AAA}\r \AAA$ and we define the
``Shukla''-complex \cite{Shukla} of $\AAA$ as
\[
\CC_{\mathrm{sh}}(\AAA,M)=\CC(\bar{\AAA},M)
\]
and as usual $\CC_{\mathrm{sh}}(\AAA)=\CC_{\mathrm{sh}}(\AAA,\AAA)$.
\begin{propositions}
\label{ref-2.4.1-16}
$\CC_{\mathrm{sh}}(\AAA,M)$ is independent of $\ovl{\AAA}$ in
   $D(k)$.
\end{propositions}
\begin{proof}
This follows from \eqref{ref-2.5-15} (applied to $\ovl{\AAA}$)
   together with Lemma \ref{ref-2.2.2-7} and Proposition \ref{ref-2.3.2-11}. The
   detailed proof is left to the reader.
\end{proof}
Proposition \ref{ref-2.4.1-16} implies that $\CC_{\mathrm{sh}}(\AAA)$ is
well-defined in $D(k)$ which is rather weak. In \S\ref{ref-4.2-30}
we will use results of Keller \cite{keller6} to explain why
$\CC_{\mathrm{sh}}(\AAA)$ is well-defined in the homotopy
category  $\Ho(B_\infty)$ of $B_\infty$-algebras and enjoys some
functoriality properties extending the ``limited functoriality'' mentioned
above.

\subsection{Non-small categories}
\label{ref-2.5-17}
The definition of $\CC_{\mathrm {sh}}(\AAA,M)$ involves products of
abelian groups which are indexed by tuples of objects of $\AAA$. This
creates a minor set theoretic problem if $\AAA$ is not $\Uscr$-small.
Therefore in such a situation we will (implicitly)  select a
larger universe $\Vscr\supset \Uscr$ such that $\AAA$ is
$\Vscr$-small. It is clear from \eqref{ref-2.4-14} that the resulting 
Hochschild
complex is, up to isomorphism, independent of the universe $\Vscr$.
\begin{remarks} In the situations we encounter below  $\CC_{\mathrm
     {sh}}(\AAA,M)$ will always have $\Uscr$-small cohomology even if
   $\AAA$ is not $\Uscr$-small (although this will not always be
   obvious). Hence $\CC_{\mathrm {sh}}(\AAA,M)$ will always be
   $\Uscr$-small in a homotopy theoretic sense.
\end{remarks}
\section{Hochschild cohomology and deformation theory of abelian
   categories}
\label{ref-3-18}
If $\UUU$ is a small linear category then
$\Mod(\UUU)=\Add(\UUU^{\op},\Mod(k))$ is the category of right
$\UUU$-modules.\footnote{In \cite{lowenvdb1} $\Mod(\UUU)$ was used to denote
the category of \emph{left} $\UUU$-modules. The category of right
$\UUU$-modules was denoted by $\Pre(\UUU)$.}

Assume that $\Ascr$ is a small $k$-linear abelian category and let
$\Mscr$ be an object in $\Mod(\Ascr^{\op}\otimes \Ascr)$.  The
category $\Cscr=\Ind \Ascr$ is the formal closure of $\Ascr$ under
filtered small colimits.  We extend $\Mscr$ to an object
$\tilde{\Mscr}$ in $\Mod(\Cscr^{\op}\otimes \Cscr)$ by
\[
\tilde{\Mscr}(\dirlim_{i\in I} A_i, \dirlim_{j\in J} B_j)=
\invlim_{i\in I} \dirlim_{j\in J} M(A_i,B_j)
\]
  It is well-known that $\Cscr$ is a Grothendieck category
and in particular it has enough injectives.
Let $\III=\Inj(\Cscr)$ be the
category of injectives of $\Cscr$. We restrict $\tilde{\Mscr}$ to an
$\III$-bimodule $M$ and we define the Hochschild complex of
$\Mscr$ as
\begin{equation}
\label{ref-3.1-19}
\CC_{\mathrm{ab}}(\Ascr,\Mscr)=\CC_{\mathrm{sh}}(\III,M)
\end{equation}
and $\CC_{\mathrm{ab}}(\Ascr)=\CC_{\mathrm{ab}}(\Ascr,\Ascr)$
(note that $\III$ is not $\Uscr$-small!)  Definition \eqref{ref-3.1-19} is
motivated by the deformation theory of abelian categories which was
introduced in \cite{lowenvdb1}. For the convenience of the reader we now sketch
this theory in order to show its relation to
\eqref{ref-3.1-19}.

Below we assume that $k$ is coherent. The coherentness assumption in
necessary in deformation theory for technical reasons (see
\cite{lowenvdb1}) but it it \emph{not} necessary for the definition 
of the Hochschild
complex \eqref{ref-3.1-19}.

To start we need a notion of flatness. We say that $\Ascr$ is
\emph{flat} if the injectives in $\Cscr$ are $k$-flat in $\Cscr^\circ$.  We
refer to \cite[\S3]{lowenvdb1} for several equivalent (but more technical)
characterizations of flatness which are intrinsic in terms of
$\Ascr$ (in particular from one of those characterizations it follows
that flatness is self dual). Below we assume that $\Ascr$ is flat.

  Let $\Rng^0$ be the category whose objects
are coherent rings and whose morphisms are surjective maps with
finitely generated nilpotent kernel and let $\theta:l\r k$ be an
object in $\Rng^0/k$.
  A \emph{flat $l$-deformation} of $\Ascr$ is
an $l$-linear flat abelian category $\Bscr$ together with an equivalence
$\Bscr_k\cong \Ascr$. Here $\Bscr_k$ is the full-subcategory of
$\Bscr$ whose objects are annihilated by $\ker \theta$.

Let $\Def_\Ascr(l)$ be the groupoid whose objects are the flat
$l$-deformations of $\Ascr$ and whose morphisms are the equivalences of
deformations (in an obvious sense) up to natural isomorphism.
The
following result shows that Hochschild cohomology describes the
obstruction theory for the functor 
\[
\Def_\Ascr:\Rng^0/k\r \Gd
\]
\begin{theorem}
\label{ref-3.1-20}
Let $\sigma:l'\r l$ be a map in $\Rng^0/k$ such that $I= \ker \sigma$
is annihilated by $\ker (l'\r k)$. Let $\Bscr\in \Def_\Ascr(l)$ and
let $\Def_\Ascr(\sigma)^{-1}(\Bscr)$ be the groupoid whose objects are
flat $l'$-deformations of $\Bscr$. 
 \begin{enumerate}
\item There is an obstruction $o(\Bscr)\in
   \mathrm{H}\CC^3_{\mathrm{ab}}(\Ascr,I\otimes \Ascr)$ such that
$
   \Ob(\Def_\Ascr(\sigma)^{-1}(\Bscr))\neq \varnothing \iff o(\Bscr)=0
$
\item If $o(\Bscr)=0$ then $\Sk(\Def_\Ascr(\sigma)^{-1}(\Bscr))$ is an
   affine space over $\mathrm{H}\CC^2_{\mathrm{ab}}(\Ascr,I\otimes \Ascr)$
\end{enumerate}
\end{theorem}
We now sketch the proof of this theorem. The
proof will be complete in the case where  $k$ is a field (and this
is already sufficient motivation for \eqref{ref-3.1-19}). In the
general case the proof depends on some results concerning the
deformation theory of DG-categories which are well-known to experts
but which do not seem to have appeared in the literature yet. We refer to
\cite{kellerggeigle,lowen2}.

Assume that $\UUU$ is a $k$-linear category. We say that $\UUU$ is
flat if all Hom-sets of $\UUU$ are $k$-flat. Assume that $\UUU$ is
flat. A flat $l$-deformation of $\UUU$ is a flat $l$-linear category
$\VVV$ together with an equivalence $\VVV\otimes_l k\cong \UUU$ where
$\VVV\otimes_l k$ is obtained by tensoring the $\Hom$-sets of $\VVV$
with $k$.   As above the morphisms between deformations are
equivalences, up to natural isomorphism. The corresponding groupoid is
denoted by $\ddef_{\UUU}(l)$.  The groupoid $\ddef^s_{\UUU}(l)$ of
\emph{strict} deformations of $\UUU$ is defined similarly except that
we replace ``equivalence'' by ``isomorphism'' everywhere.

We recall the following:
\begin{proposition} {\cite[Thm B.4]{lowenvdb1}} \label{ref-3.2-21} The natural functor 
  $\ddef^s_{\UUU}(l)\r \ddef_{\UUU}(l)$ defines a bijection between
  the corresponding skeletons.
\end{proposition}

As above let $\III$ be the category of injectives in $\Cscr$. The
following is one of the main results of \cite{lowenvdb1}. 
\begin{theorem} \label{ref-3.3-22} \cite[Thm.\ 8.8,Thm.\ 
8.17]{lowenvdb1} The category $\III$ is flat as linear
   category and there is an equivalence of categories between 
   $\Def_\Ascr(l)$ and $\ddef_{\III}(l)$.
\end{theorem}
It follows that the deformation theory of abelian categories
reduces to the deformation theory of linear categories.

Let $\UUU$ be a flat $k$-linear category as above.  If $k$ is a field
then it is well-known that (as in the algebra case) the strict
deformation theory of $\UUU$ is controlled by the Hochschild
cohomology of $\UUU$. Let $\sigma:l'\r l$ be a map in $\Rng^0/k$ such
that $I= \ker \sigma$ is annihilated by $\ker (l'\r k)$ and let $\VVV$
be a flat $l$-deformation of $\UUU$.  Then it follows from Proposition
\ref{ref-3.2-21} (replacing $k$ by $l$ and $l$ by $l'$) that there is a bijection
$\Sk(\ddef^s_{\UUU}(\sigma)^{-1}(\VVV))\cong
\Sk(\ddef_{\UUU}(\sigma)^{-1}(\VVV))$. Hence the non-strict
deformation theory of $\UUU$ is 
controlled by Hochschild cohomology as well and this then leads to a
proof of \eqref{ref-3.1-20}.

If $k$ is not a field then there
is a technical difficulty in the sense that when the $\Hom$-sets of
$\UUU$ are not projective over $k$, the naive Hochschild cohomology of
$\UUU$, $\mathrm{H}\CC^\ast(\UUU,I \otimes \UUU)$ does \emph{not} lead
to the correct obstruction theory for the deformations of $\UUU$. This
problem is rather serious since  the linear category
$\III=\Inj(\Cscr)$ will have flat, but not in general projective
$\Hom$-sets.

Nevertheless it is true that the deformation theory of $\UUU$ is
controled by $\mathrm{H}\CC_{\mathrm{sh}}^\ast(\UUU)$ \cite{lowen2}.
To prove this one replaces $\UUU$ by an  appropriate semi-free resolution \cite[Lemma 13.5]{Drinfeld2}
$\bar{\UUU}$ and one studies the deformations of $\bar{\UUU}$ in the
homotopy category of DG-categories which turns out to be controlled by
the homology of $\Der(\bar{\UUU},I\otimes \bar{\UUU})$ in degree 1 and
2. There is an exact triangle 
\[
\Der(\ovl{\UUU},I\otimes\bar{\UUU})[-1]\r \CC(\ovl{\UUU},I\otimes \ovl{\UUU})\r
 I\otimes \ovl{\UUU} \r
\]
and therefore in degrees $\ge 2$, $\CC(\ovl{\UUU},I\otimes\ovl{\UUU})$ and
$\Der(\ovl{\UUU},I\otimes \bar{\UUU})[-1]$ have the same  homology which
finished the proof.
\begin{remark} If $\UUU$ is flat and $M$ is a $\UUU$ bimodule then it
     is not hard to see that there is a quasi-isomorphism in $D(k)$
\[
\CC_{\mathrm{sh}}(\UUU,M)=\RHom_{\UUU^{\op}\otimes \UUU}(\UUU,M)
\]
This ``more elementary'' interpretation of  Shukla cohomology does not
seem to be useful for deformation theory however.
\end{remark}
\begin{remark}
   If $k$ is a field of characteristic zero then it follows from 
   Proposition \ref{ref-3.2-21} and 
   Theorem \ref{ref-3.3-22} that there is a bijection between
   $\Sk(\Def_\Ascr(l))$ and the solutions of the Maurer-Cartan equation
   in $\CC_{\mathrm{ab}}(\Ascr)[1]$ with coefficients in $\ker(l\r k)$, modulo
   gauge equivalence \cite{Ko3}. Hence $\CC_{\mathrm{ab}}(\Ascr)[1]$ is
   the
DG-Lie algebra controling the deformation theory of $\Ascr$.
\end{remark}
\begin{remark}  Assume that $k$ is a field. If
  $\dim \mathrm{H}\CC^2_{\mathrm{ab}}(\Ascr)< \infty$ then $\Ascr$ has
  a (formal) versal deformation \cite[Chapt.\ 
  6]{versal,Stevens}. Or in equivalent terms: the functor
\[
\Sk(\Def_\Ascr(-)):\Rng^0/k\r \Set
\]has a hull \cite{Schlessinger}. I.e.\ there is
  a noetherian local ring $(R,m)$ with residue field $k$ together with
  compatible flat deformations $\Ascr_n$ of $\Ascr$ over $R/m^n$ such
  that the formal object $\invlim_n \Ascr_n$ in $\Def_\Ascr(R)$
  satisfies the versality condition. Using the existence criterion in
  \cite{Schlessinger} this can be shown by translation to the linear case
  using Proposition \ref{ref-3.2-21} and Theorem \ref{ref-3.3-22}.

Of course one actually wants an $R$-linear abelian category $\tilde{\Ascr}$
representing the formal object $\invlim_n \Ascr_n$. This is possible
under certain conditions. See \cite{lowen3}.
\end{remark}
\begin{remark}
In the rest of this paper we restrict ourselves to the study of
$\CC_{\mathrm{ab}}(\Ascr)$ as this is the most interesting case. And in any
case this is sufficient for deformation theory in case $k$ is a field.
\end{remark}
\section{More background on Hochschild cohomology of
   DG-categories}
\label{ref-4-23}
As we have seen the Hochschild cohomology of an
abelian category is basically the Hochschild cohomology of a suitable
DG-category. So we need techniques for computing the Hochschild
cohomology of DG-categories. A powerful tool in this respect was
provided by Keller in \cite{keller6}.
\subsection{Keller's results}
\label{ref-4.1-24}
All results in this section are due to Keller. See \cite{keller6}.
Suppose we have small $k$-cofibrant DG-categories $\AAA$ and $\BBB$ and
suppose that $X$ is a cofibrant  $\AAA^{\op} \otimes \BBB$-module. Then by
functoriality, we get a map of $\AAA{-}\AAA$-bimodules
\[
\lambda: \AAA
\lra \Hom_{\BBB}(X,X) = \RHom_{\BBB}(X,X)
\]
  and an induced map
\[
\lambda_{\ast}: \CC(\AAA) \lra \CC(\AAA,\Hom_{\BBB}(X,X))
\]
Similarly,
\[
\omega: \BBB \lra \Hom_{\AAA^{\op}}(X,X) = \RHom_{\AAA^{\op}}(X,X)
\] induces
\[
\omega_{\ast}: \CC(\BBB) \lra \CC(\BBB,\Hom_{\AAA^{\op}}(X,X))
\]
Below, in case $X$ is variable, we will adorn $\lambda$ and
$\omega$ by a subscript $X$.

Let $\CCC$ be the
DG-category such that
\[
\Ob(\CCC)=\Ob(\AAA)\coprod \Ob(\BBB)
\]
and such that
\[
\CCC(U,V)=
\begin{cases}
\AAA(U,V) &\text{if $U,V\in \AAA$}\\
\BBB(U,V) &\text{if $U,V\in \BBB$}\\
X(U,V)&\text{if $U\in \BBB$, $V\in \AAA$}\\
0&\text{otherwise}
\end{cases}
\]
Below we will sometimes use the notation
$(\AAA\xleftarrow{X}\BBB)$ for the category $\CCC$ and we refer to
$\CCC$
as an ``arrow category''.

Consider the canonical inclusions $i_{\AAA}: \AAA \lra \CCC$ and
$i_{\BBB}: \BBB \lra \CCC$.
\begin{theorems}\label{ref-4.1.1-25}\cite[\S 4.4,4.5]{keller6}
\begin{enumerate}
\item There is a quasi-isomorphism $\CC(\AAA,\RHom_{\BBB}(X,X)) \cong
\CC(\BBB,\RHom_{\AAA^{\op}}(X,X))$.
\item If $\lambda_{\ast}$ is a quasi-isomorphism, $i_{\BBB}^{\ast}: \CC(\CCC)
\lra \CC(\BBB)$ is a quasi-isomorphism of $B_{\infty}$-algebra's.
\item If $\omega_{\ast}$ is a quasi-isomorphism, $i_{\AAA}^{\ast}: \CC(\CCC)
\lra \CC(\AAA)$ is a quasi-isomorphism of $B_{\infty}$-algebra's.
\end{enumerate}
\end{theorems}
So in particular if  $\lambda_\ast$ is a  quasi-isomorphism then we have an
induced map in $\Ho(B_\infty)$
\[
\phi_X:\CC(\BBB)\r \CC(\AAA)
\]
given by $i_\AAA^\ast (i_\BBB^\ast)^{-1}$. If  $\omega_\ast$ is a 
   quasi-isomorphism then $\phi_X$ is an isomorphism. 
\begin{theorems}\cite{keller6}\label{ref-4.1.2-26}
\begin{enumerate}
\item $\phi_X$ depends only on the isomorphism class of $X$ in
   $D(\AAA^{\op}\otimes \BBB)$.
\item If $j:\AAA\r \BBB$ is fully faithful and $X$ is given by
   $X(B,A)=\AAA(B,j(A))$ then $\lambda$ is a quasi-isomorphism and
   $\phi_X=j^\ast$ in $\Ho(B_\infty)$.
\end{enumerate}
\end{theorems}
\begin{proof} 1.\ is proved as \cite[Thm. 4.6a]{keller6}. 2.\ is
  \cite[Thm. 4.6c]{keller6}.
\end{proof}
\begin{theorems} \cite[Thm 4.6b]{keller6} \label{ref-4.1-27} If 
$-\otimes_\AAA X$ induces a fully
   faithful functor $D(\AAA)\r D(\BBB)$ then $\lambda$ is a
   quasi-isomorphism and hence
there is a well defined
   associated ``restriction'' morphism $\phi_X:\CC(\BBB)\r
   \CC(\AAA)$. If $-\otimes_\AAA X$ induces an equivalence then
   $\phi_X$ is an isomorphism.
\end{theorems}
There is also a transitivity result for the maps $\phi_X$.
\begin{theorems} \cite{keller6} \label{ref-4.1.3-28} Let $X$ be a
   cofibrant $\AAA^{\op}\otimes \BBB$-module
   such that $(\lambda_X)_\ast$ is a  quasi-isomorphism and let $Y$ be a cofibrant
   $\BBB^{\op}\otimes \CCC$-module inducing a fully faithful
   functor $D(\BBB)\r D(\CCC)$. Then $(\lambda_{X\otimes_\BBB Y})_\ast$ is an
   isomorphism and $\phi_{X\otimes_\BBB Y} =\phi_Y\circ \phi_X$.
\end{theorems}
\begin{proof} This is proved as \cite[Thm 4.6d]{keller6}.
\end{proof}
\begin{remarks}
\label{ref-4.1.4-29} If $j:\AAA\r \BBB$ is fully faithful then instead of
   the $\AAA-\BBB$-module $X$ defined in Theorem \ref{ref-4.1.2-26} it is
   just as natural to use the $\BBB-\AAA$ bimodule $Y$ defined by
   $Y(A,B)=\Hom_\BBB(j(A),B)$.  One may dualize all arguments for this
   bimodule.  The relevant arrow category is now
   $\CCC'=(\BBB\xleftarrow{Y} \AAA)$ where the arrow is just the
   inclusion. It turns our that now $\omega_Y$
   is a  quasi-isomorphism. Thus we obtain a morphism in $\Ho(B_\infty)$
given
   by $i_\BBB^\ast (i_\AAA^{\ast})^{-1}$. Dualizing the proof of
   Theorem \ref{ref-4.1.2-26}.2 above we still find $i_\BBB^\ast
   (i_\AAA^{\ast})^{-1}=j^\ast$ in $\Ho(B_\infty)$. Hence the ``bimodule
   interpretation'' of limited functoriality is unambigous.
\end{remarks}
\subsection{The functoriality of the
   Shukla complex}
\label{ref-4.2-30}
Now we do not assume that our DG-categories are $k$-cofibrant.
Using the results in \S\ref{ref-4.1-24} we can explain why the Shukla
complex is well-defined and functorial.  More precisely we show that
$\CC_{\mathrm{sh}}(-):\AAA\r \CC_{\mathrm{sh}}(\AAA)$ defines a
contravariant functor on a suitable category of  small DG-categories
with values in $\Ho(B_\infty)$.

Let $\Fscr$ be the category whose objects are small DG-categories. If
$\AAA,\BBB\in \Ob(\Fscr)$ then $\Fscr(\AAA,\BBB)$ is defined as the set of
equivalence classes of triples $(\bar{\AAA},X,\bar{\BBB})$ where
$\bar{\AAA}\r \AAA$, $\bar{\BBB}\r \BBB$ are $k$-cofibrant resolutions
and $X$ is a cofibrant $\bar{\AAA}^{\op}\otimes \bar{\BBB}$-module
such that $-\otimes_{\bar{\AAA}} X$ induces a fully faithful functor
$D(\AAA)\cong D(\bar{\AAA})\r D(\bar{\BBB})\cong D(\BBB)$. Two triples
$(\bar{\AAA},X,\bar{\BBB})$, $(\bar{\AAA}',X',\bar{\BBB}')$ are
equivalent if $X$ and $X'$ correspond under the canonical equivalence
between $D(\bar{\AAA}^{\op}\otimes \bar{\BBB})$ and
$D(\bar{\AAA}^{\prime \op}\otimes \bar{\BBB}')$.

If $\bar{\AAA}\r \AAA$, $\bar{\AAA}'\r \AAA$ are $k$-cofibrant
resolutions then we define $C_{\bar{\AAA}\bar{\AAA}'}$ to be a
cofibrant $\bar{\AAA}^{\op}\otimes_k \bar{\AAA}^{\prime}$ resolution
of $\AAA$ considered as
$\bar{\AAA}^{\op}-\bar{\AAA}^{\prime}$-bimodule.
$(\bar{\AAA},C_{\bar{\AAA}\bar{\AAA}'},\bar{\AAA}')$ defines a 
canonical equivalence
class of objects in
$\Fscr(\AAA,\AAA)$ which we denote by $\Id_\AAA$ since
$C_{\bar{\AAA}\bar{\AAA}'}$ induces the identity on $D(\AAA)$.

Composition of
triples  is defined as
\[
(\bar{\AAA},X,\bar{\BBB})\circ
(\bar{\BBB}',Y,\bar{\CCC})=(\bar{\AAA},X\otimes_{\bar{\BBB}}
C_{\bar{\BBB}\bar{\BBB}'}\otimes_{\bar{\BBB}'} Y,\bar{\CCC})
\]
It is easy to see that this is compatible with
equivalence. Furthermore the maps $\Id_\XXX$ for $\XXX\in \Ob(\Fscr)$
behave as identities for this composition.

Now for every $\AAA\in \Ob(\Fscr)$ fix a $k$-cofibrant resolution
$\bar{\AAA}$ and define $\CC_{\mathrm{sh}}(\AAA)=\CC(\bar{\AAA})$. We
will make $\CC_{\mathrm{sh}}(-)$ into a functor on $\Fscr$.  Assume we
have a triple $(\bar{\AAA}',X',\bar{\BBB}')$. Then we have an
equivalent triple of the form $(\bar{\AAA},X,\bar{\BBB})$, with
$X=C_{\bar{\AAA}\bar{\AAA}'}\otimes_{\bar{\AAA}'}X'
\otimes_{\bar{\BBB}'} C_{\bar{\BBB}'\bar{\BBB}}$.  To the triple
$(\bar{\AAA}',X',\bar{\BBB}')$ we now associate the map
$\phi_X:\CC(\bar{\AAA})\r \CC(\bar{\BBB})$. Using Theorems
\ref{ref-4.1.2-26},\ref{ref-4.1.3-28} it is easy to see that the assignment
$(\bar{\AAA}',X',\bar{\BBB}')\mapsto \phi_X$ is compatible with
equivalence, compositions and sends 
 $\Id_\AAA$ to the identity.
In this way we have reached our goal of making $\CC_{\mathrm{sh}}(-)$
into a functor.
\begin{remarks} The functor $\CC_{\mathrm{sh}}(-)$ inverts
   quasi-equivalences.  This follows from the corresponding result in
   the $k$-cofibrant case (see  Theorem  \ref{ref-4.1-27}).
\end{remarks}
\begin{remarks}
If we
choose another system of resolutions $q:\bar{\AAA}'\r \AAA$ then
we have a canonical isomorphism $\CC(\bar{\AAA})\r
\CC(\bar{\AAA}')$ induced by $C_{\bar{\AAA}\bar{\AAA}'}$
which defines a natural isomorphism between the
functors  $\AAA\longmapsto \CC(\bar{\AAA})$ and  $\AAA\longmapsto
\CC(\bar{\AAA}')$. So the functor
$\CC_{\mathrm{sh}}(-)$ is well defined  up  to a canonical natural isomorphism.
\end{remarks}
\begin{remarks}
\label{ref-4.2.3-31}
If $j:\AAA\r \BBB$ is a fully faithful functor and if $\bar{\BBB}\r
\BBB$ is a $k$-cofibrant resolution of $\BBB$ then we may restrict
this resolution to a $k$-cofibrant resolution of $\AAA$. So $j$
extends to a fully faithful functor $\bar{j}:\bar{\AAA}\r
\bar{\BBB}$. Thus $\bar{j}$ defines a morphism of Shukla complexes
\[
\CC_{\mathrm{sh}}(\BBB)\r\CC_{\mathrm{sh}}(\AAA)
\]
which we will denote by $j^\ast$. This notation is natural by Remark
\ref{ref-4.1.4-29}.
\end{remarks}
\subsection{The ``Cosmic Censorship'' principle}
\label{ref-4.3-32}
Assume that $\AAA$ and $\BBB$ are small $k$-cofibrant DG-categories.
It follows from Theorem \ref{ref-4.3.3-36} that $\CC(\AAA) \cong
\CC(\BBB)$ in $\Ho(B_{\infty})$ if $\omega$ and $\lambda$ are
quasi-isomorphims (and indeed this is the way the result is explicitly
stated in \cite{keller6}). $\lambda$ and $\omega$ being
quasi-isomorphisms is equivalent to
\begin{align*}
\lambda(A,A'): \AAA(A,A')&\r \RHom_\BBB(X(-,A),X(-,A'))\\
\omega(B,B'):\BBB(B,B')&\r \RHom_{\AAA^{\op}}(X(B',-),X(B,-))
\end{align*}
being quasi-isomorphisms for all $A,A'\in \Ob(\AAA)$ and $B,B'\in
\Ob(\BBB)$.

The extra generality of Theorem
\ref{ref-4.1.1-25} will be essential for us when we study ringed spaces,
for it turns out that sometimes $\lambda_\ast$ or $\omega_\ast$ are
quasi-isomorphisms when this is not necessarily the case for $\lambda$
or $\omega$.

In the application to ringed spaces $\Ob(\AAA)$ will be equipped with
a non-trivial transitive relation $\rrr$ such that $\AAA(A,A') = 0$ if $(A,A')
\notin \rrr$. We will call such $\Rscr$ a \emph{censoring relation}.
Note that any $\AAA$ has a \emph{trivial censoring relation} given by
$\Rscr=\Ob(\AAA)\times \Ob(\AAA)$.
We
have the following result.
\begin{propositions}
\label{ref-4.3.1-33}
Assume that $\AAA$ has a censoring relation $\Rscr$ and let $M$
be an $\AAA{-}\AAA$-bimodule. Define $M_0$ by
\[
M_0(A,A')=
\begin{cases}
M(A,A')&\text{if $(A,A')\in \Rscr$}\\
0&\text{otherwise}
\end{cases}
\]
Then $M_0$ is a subbimodule of $M$ and
furthermore
\begin{equation}
\label{ref-4.1-34}
\CC(\AAA,M)=\CC(\AAA,M_0)
\end{equation}
\end{propositions}
\begin{proof}
That $M_0$ is a submodule is clear. The equality
$\CC(\AAA,M)=\CC(\AAA,M_0)$ follows immediately from the definition of
the Hochschild complex.
\end{proof}
\begin{propositions}
\label{ref-4.3.2-35}
Assume that $\AAA$ has a censoring relation $\Rscr$ and
\[
\lambda(A,A'): \AAA(A,A')\r \RHom_\BBB(X(-,A),X(-,A'))
\]
is an isomorphism for all $(A,A')\in \Rscr$. Then $\lambda_\ast$
is a quasi-isomorphism.
\end{propositions}
So in a sense the zero $\Hom$-sets in $\AAA$ censor  possible ``bad parts'' of
$\RHom_\BBB(X,X)$. See Remark \ref{ref-7.3.2-84} below for an application of
this principle.
\begin{proof}
Our hypotheses imply that  $\lambda$ factors as
\[
  \AAA(A,A')\r \Hom_\BBB(X,X)_0\r \Hom_\BBB(X,X)
\]
and the first map is quasi-isomorphism. Hence $\lambda_\ast$ factors
as a composition of a quasi-isomorphism and an isomorphism
\[
\CC(\AAA(A,A'))\r \CC(\AAA,\Hom_\BBB(X,X)_0)\r \CC(\AAA,\Hom_\BBB(X,X))
\]
finishing the proof.
\end{proof}

\medskip

Now assume that $\AAA$ and $\BBB$ be are arbitrary small (not
necessarily $k$-cofibrant) DG-categories which are equipped with
(possibly trivial) censoring relations $\Rscr$ and $\Lscr$ and let $X$
be an object in $\Dif(\AAA^{\op}\otimes \BBB)$. We will use the
following criterion to compare the Shukla complexes of $\AAA$ and
$\BBB$.  
\begin{propositions}
\label{ref-4.3.3-36}
Assume that the compositions 
\begin{equation}
\label{ref-4.2-37}
\begin{array}{ccccc}
\scriptstyle  \Lambda_{A,A'}: \AAA(A,A')&\scriptstyle \r& \scriptstyle \Hom_\BBB(X(-,A),X(-,A')) 
&\scriptstyle \mathbin{\mathop{\longrightarrow}\limits^{\mathrm{can}}} 
&\scriptstyle \RHom_\BBB(X(-,A),X(-,A'))\\
 \scriptstyle   \Omega_{B,B'}: \BBB(B,B')&\scriptstyle \r& \scriptstyle \Hom_{\AAA^{\op}}(X(B',-),X(B,-))
& \scriptstyle \mathbin{\mathop{\longrightarrow}\limits^{\mathrm{can}}}&
\scriptstyle    \RHom_{\AAA^{\op}}(X(B',-),X(B,-))
\end{array}
\end{equation}
are quasi-isomorphisms for all $(A,A')\in \Rscr$ and $(B,B')\in
\Lscr$. Then $\CC_{\mathrm{sh}}(\AAA)\cong \CC_{\mathrm{sh}}(\BBB)$ in
$\Ho(B_\infty)$.
\end{propositions}
\begin{proof}
  Let $\bar{\AAA} \lra
\AAA$ and $\bar{\BBB} \lra \BBB$ be $k$-cofibrant resolutions. By
replacing
$\bar{\AAA}$ and $\bar{\BBB}$ by $\bar{\AAA}_0$ and $\bar{\BBB}_0$
respectively we may \emph{and we will} assume that $\Rscr$ and $\Lscr$
are censoring relations for $\bar{\AAA}$ and $\bar{\BBB}$.

Now let $\bar{X}
\lra X$ be a cofibrant resolution of $X$ in $\Dif({\ovl{\AAA}}^{\op}
\otimes \ovl{\BBB})$. We may apply Theorem \ref{ref-4.1.1-25}
and Proposition \ref{ref-4.3.2-35}
  to
the triple $(\ovl{\AAA},\ovl{X},\ovl{\BBB})$. Thus we need to check
that the appropriate $\lambda(A,A')$ and $\omega(B,B')$ are
quasi-isomorphisms.

For $A,A'\in \Ob(\ovl{\AAA})$ we have a commutative diagram in $D(k)$
\[
\xymatrix{\ovl{\AAA}(A,A') \ar[r]^-{\lambda(A,A')} \ar[d]_{\cong} &
{\Hom_{\ovl{\BBB}}(\bar{X}(-,A),\bar{X}(-,A'))} \ar@{.}[d]^{\cong}\\
{\AAA(A,A')} \ar[r]_-{\Lambda_{A,A'}} &{\RHom_{\BBB}(X(-,A),
   X(-,A'))}}
\]
in which the quasi-isomorphism to the right is justified by Lemma
\ref{ref-2.2.2-7}. Hence $\lambda(A,A')$ is a quasi-isomorphism if this is
true for the corresponding map in \eqref{ref-4.2-37}. A similar
observation holds for $\omega(B,B')$. This finishes the proof.
\end{proof}
We sometimes use the following compact criterion to see if a fully
faithful map induces an isomorphism on Shukla cohomology.
\begin{propositions}
\label{ref-4.3.4-38} Let $j:\AAA\r \BBB$ be fully faithful and  assume that
$\BBB$ has a (possibly trivial) censoring relation $\Lscr$ such that
for all $(B,B')\in \Rscr$ the canonical maps given by functoriality
\[
\BBB(B,B')\r \RHom_{\AAA^{\op}}(\Hom_\BBB(B',j(-)),\Hom_\BBB(B,j(-)))
\]
are isomorphisms then $j^\ast:\CC_{\mathrm{sh}}(\AAA) \r 
\CC_{\mathrm{sh}}(\BBB)$ is an
isomorphism in $\Ho(B_\infty)$.
\end{propositions}
\begin{proof} As in Remark \ref{ref-4.2.3-31} we may lift 
$j:\AAA\r\BBB$ to a fully
   faithfully map $\bar{j}:\bar{\AAA}\r \bar{\BBB}$ between
   $k$-cofibrant resolutions of $\AAA$ and $\BBB$.  Let $\bar{X}$ be a 
$\bar{\AAA}-\bar{\BBB}$-bimodule
   which is a cofibrant resolution of the bimodule given by
   $X(B,A)=\Hom_\BBB(A, j(B))$. Then by Theorem \ref{ref-4.1.2-26}
   $\phi_{\bar{X}}=\bar{j}^\ast$ and $\lambda$ is a quasi-isomorphism.
   So we have to check that $\omega_\ast$ is a quasi-isomorphism. As in
   the proof of Proposition \ref{ref-4.3.3-36} we may do this by checking
   that the $\Omega_{B,B'}$ are quasi-isomorphisms for $(B,B')\in
   \Lscr$. But this is precisely the assertion of the current
   proposition.
\end{proof}
In other words for the previous proposition to apply the contravariant
representable functors corresponding to objects in $\BBB$ should have
the same $\RHom$'s as their restrictions to $\AAA$. Note that
Proposition \ref{ref-4.3.4-38} has a dual version using contravariant
representable functors which we will also use below.
\theoremstyle{plain} \newtheorem*{convention}{Convention}
\begin{convention} Below, unless otherwise specified, if we write
$\CC_1\cong \CC_2$ for $B_{\infty}$-algebras $\CC_1$ and $\CC_2$
we mean that $\CC_1$ and $\CC_2$ are isomorphic in $\Ho(B_{\infty})$.
\end{convention}
\subsection{DG-categories of cofibrant objects}
In this section we give an easy but rather spectacular application of
Proposition \ref{ref-4.3.4-38}. A weak version of it will be used
afterwards to compare the Hochschild cohomology of an abelian category
and its bounded derived category (viewed as a DG-category).

Let $\AAA$ be a small DG-category and let
\[
\AAA\r \Dif(\AAA):\AAA\r \AAA(-,A)
\]
be the Yoneda embedding. We have the following
\begin{theorems}
\label{ref-4.4.1-39}
Let $\BBB$ be any DG-subcategory of $\Dif(\AAA)$
   which consists of cofibrant objects and which contains $\AAA$.
Then the restrictiom map $\CC_{\mathrm{sh}}(\BBB)\r
\CC_{\mathrm{sh}}(\AAA)$ is a quasi-isomorphism.
\end{theorems}
\begin{proof}
   Let $j:\AAA\r \BBB$ be the inclusion functor.  By the dual version
   of Proposition \ref{ref-4.3.4-38} we need to prove
\begin{equation}
\label{ref-4.3-40}
\RHom_{\AAA}(\BBB(j(-),B), \BBB(j(-),B'))= \BBB(B,B')=\Hom_{\AAA}(B,B')
\end{equation}
  for
$B,B'\in \BBB$. Since $\BBB(j(-),B)=\Hom_\AAA(j(-),B)=B$ for $B\in
   \BBB$ (recall $\Dif(\AAA)=\DGFun(\AAA^{\op},C(k))$),
  \eqref{ref-4.3-40} follows from the fact that $B$ is cofibrant.
\end{proof}
\subsection{Hochschild cohomology as a derived center}
\label{ref-4.5-41}
It is well known that  the Hochschild cohomology of a ring may be regarded
as a kind of non-additive derived version of the center. In this
section we show that this generalizes trivially to Shukla cohomology
of DG-categories. That is, we will construct for a small DG-category
$\AAA$ a canonically defined map
\[
\sigma_\AAA:H\CC_{\mathrm{sh}}^\ast(\AAA)\r Z(H^\ast(\AAA))
\]
where $Z(H^\ast(\AAA))$ is the center of the graded category
$H^\ast(\AAA)$.

Assume that $\UUU$ is a small $k$-linear $\ZZ$-graded category. The
\emph{center} $Z(\UUU)$ of $\UUU$ is by definition the ring of graded
endomorphisms of the identity functor on $\UUU$. More concretely the
center of $\UUU$ is a graded ring whose homogeneous elements consists
of tuples of homogeneous elements $(\phi_U)_U\in \prod_{U\in
   \UUU}\UUU(U,U)$ such that for any homogeneous $f\in \UUU(U,V)$ one
has $f\phi_U=(-1)^{|f||\phi_U|} \phi_V f$. In particular $Z(\UUU)$ is
(super) commutative.

Assume that $\AAA$ is a $k$-cofibrant DG-category.  Let
\[
\Sigma_\AAA:\CC(\AAA)\r \prod_{A\in \AAA}\AAA(A,A)
\]
be the map associated to the morphism of double complexes which sends
$\DD(\AAA)$ to its first column. An easy computation with the explicit
formulas for the cup product \cite{getzlerjones} shows that
$\sigma_\AAA\overset{\text{def}}{=}H^\ast(\Sigma_\AAA)$ is a graded ring map
which maps $H\CC^\ast(\AAA)$ to the center of $H^\ast(\AAA)$

Now let $\AAA$ be an arbitrary small DG-category and let $\bar{\AAA}\r
\AAA$ be a $k$-cofibrant resolution of $\AAA$. We define
$\bar{\sigma}_{\bar{\AAA}}:H\CC^\ast_{\mathrm{sh}}(\AAA) \r
Z(H^\ast(\AAA))$ as the composition
\[
\CC_{\mathrm{sh}}(\AAA)=
\CC(\bar{\AAA})\xrightarrow{\sigma_{\bar{\AAA}}}
Z(H^\ast(\bar{\AAA}))\cong  Z(H^\ast(\AAA)))
\]
Let $\bar{\AAA}'\r \AAA$ be another
$k$-cofibrant resolution of $\AAA$ and let $C_{\bar{\AAA}\bar{\AAA}'}$
be as in \S\ref{ref-4.2-30}.  Let
$\bar{\CCC}=(\bar{\AAA}\xleftarrow{C_{\bar{\AAA}\bar{\AAA}'}}\bar{\AAA}')$.
Then
\[
H^\ast(\bar{\CCC})=(H^\ast(\bar{\AAA})
\xleftarrow{H^\ast(C_{\bar{\AAA}\bar{\AAA}'})}H^\ast(
\bar{\AAA}'))
\]
Now by construction we have canonical compatible isomorphisms
$H^\ast(\bar{\AAA})\cong H^\ast(\AAA)$, $H^\ast(\bar{\AAA}')\cong
H^\ast(\AAA)$,  and
$H^\ast(C_{\bar{\AAA}\bar{\AAA}'})\cong H^\ast(\AAA)$. Put $\CCC=(\AAA
\xleftarrow{\AAA}\AAA)$.

We then have the following commutative
diagram
\[
\begin{CD} Z(H^\ast(\AAA))@<<< Z(H^\ast(\CCC)) @>>>
Z(H^\ast(\AAA))\\
@| @| @|\\
Z(H^\ast(\bar{\AAA}))@<<< Z(H^\ast(\bar{\CCC})) @>>>
Z(H^\ast(\bar{\AAA}'))\\
@A\sigma_{\bar{\AAA}} AA @A\sigma_{\bar{\CCC}} AA
@A\sigma_{\bar{\AAA}'}AA\\
H\CC^\ast(\bar{\AAA}) @<<\cong < H\CC^\ast(\bar{\CCC}) @>>\cong> H\CC^\ast(\bar{\AAA}')
\end{CD}
\]
Using the definition of $\CCC$ it is trivial to see that the topmost
horizontal arrows are isomorphisms and compose to the identity
isomorphism
$Z(H^\ast(\AAA))\r Z(H^\ast(\AAA))$. In this way we obtain a
commutative diagram 
\begin{equation}
\label{ref-4.4-42}
\begin{CD}
Z(H^\ast(\AAA)) @= Z(H^\ast(\AAA))\\
@A\bar{\sigma}_{\AAA}  AA @AA\bar{\sigma}_{\AAA'} A \\
H\CC^\ast(\bar{\AAA}) @>>\cong> H\CC^\ast(\bar{\AAA}')
\end{CD}
\end{equation}
where the bottom horizontal arrow is coming from the canonical
$\Ho(B_\infty)$ isomorphism between $\CC^\ast(\bar{\AAA})$ and
$\CC^\ast(\bar{\AAA}')$. We now put $\sigma_\AAA=\sigma_{\bar{\AAA}}$.
Diagram \eqref{ref-4.4-42} shows that $\sigma_\AAA$ is indeed
well-defined in the appropriate sense.
\section{Grothendieck categories}\label{ref-5-43}
It follows from the definition \eqref{ref-3.1-19} that we need to be able
to understand the Hochschild cohomology of the category of injectives
in a Grothendieck category. In this section we will prove the relevant
technical results.
\subsection{A model structure}
\label{ref-5.1-44}
We assume troughout that $\ccc$ is a $k$-linear Grothendieck category
and as usual $C(\Cscr)$ denotes the category of complexes over
$\Cscr$. In this section we construct a generalization of the usual
injective model structure on $C(\Cscr)$ \cite{ALS,beke,Franke}. It is used
for some of the proofs below, but \emph{not} for the statement of the
results.

  For a small
$k$-DG-category
$\AAA$ consider the
category
$\Dif(\AAA,\Cscr)=\DGFun(\AAA^{\op},C(\ccc))$.
\begin{proposition}\label{ref-5.1-45} \begin{enumerate} \item
$\Dif(\AAA,\Cscr)$ has the structure of a model category in which the
weak equivalences are the pointwise quasi-isomorphisms and the
cofibrations are the pointwise monomorphisms.
\item Suppose
$\AAA$ is $k$-cofibrant. If
$F
\in
\Dif(\AAA,\Cscr)$ is fibrant, then so is every
$F(A)$ in $C(\ccc)$.
\end{enumerate}
\end{proposition}
If $\AAA=k$ then we obtain the usual injective model structure on
$C(\Cscr)$. A very efficient proof for the existence of the latter
has been given in \cite{beke} by Beke. Our proof of Proposition 
\ref{ref-5.1-45}
is based on the following result which is an abstraction of the method
used for $C(\Cscr)$ by Beke.
\begin{proposition}\label{ref-5.2-46} \cite{beke}
Let $(H^i)_{i \in \Z}: \aaa \lra \bbb$ be additive functors between
Grothendieck categories such that
\begin{enumerate}
\item $H^i$ preserves filtered colimits;
\item $(H^i)_i$ is
effaceable (i.e. every $A \in \aaa$ admits a mono $m: A \lra A'$ with $H^i(m) =
0$ for every $i \in \Z$, or, equivalently, $H^i(\Inj(\aaa)) = 0$ for every $i
\in \Z$);
\item $(H^i)_i$ is cohomological (i.e. for every short exact $0 \lra A' \lra A
\lra A''
\lra 0$ is $\aaa$ there is a long exact $\dots \lra H^i(A') \lra H^i(A) \lra
H^i(A'') \lra H^{i+1}(A') \lra \dots$ in $\bbb$)
\end{enumerate}
Let
$\mathsf{iso}_{\bbb}$ be the class of isomorphisms in $\bbb$ and
$\mathsf{mono}_{\aaa}$ the class of monomorphisms in $\aaa$. There is a model
structure on
$\aaa$ such that
\begin{enumerate}
\item $\mathsf{mono}_{\aaa}$ is the class of cofibrations;
\item $\cap_{i \in \Z}(H^i)^{-1}(\mathsf{iso}_{\bbb})$ is the class of weak
equivalences.
\end{enumerate}
\begin{proof} Our notations in this proof are the ones used in \cite{beke}.
We use \cite[Theorem 1.7]{beke} which is attributed
to Jeffrey Smith. By \cite[Proposition 1.12]{beke}, $\mathsf{mono}_{\aaa} =
\mathsf{cof}(I)$ for some set $I \subset \mathsf{mono}_{\aaa}$. We check the
hypotheses for Theorem \cite[Theorem 1.7]{beke}  with $W =
\cap_{i \in \Z}(H^i)^{-1}(\mathsf{iso}_{\bbb})$. (c0) and (c3) are automatic
(using \cite[Proposition 1.18]{beke}). For (c2), we are to show that
$\mathsf{inj}(\mathsf{mono}_{\aaa})
\subset W$. So consider $f: X \lra Y$ in $\mathsf{inj}(\mathsf{mono}_{\aaa})$.
It is easily seen using the lifting property of $f$ that $f$ is a split
epimorphism with an injective kernel $K$. The result follows from the long
exact sequence associated to $0 \lra K \lra X \lra Y \lra 0$ in which each
$H^i(K) = 0$. For (c2), $\mathsf{mono}_{\aaa} \cap \cap_{i \in
\Z}(H^i)^{-1}(\mathsf{iso}_{\bbb})$ is easily seen to be closed under pushouts
using the long exact sequence, and under transfinite composition using that
filtered colimits are exact in $\aaa$ and that $H^i$ preserves them.
\end{proof}
\end{proposition}
\begin{proof}[Proof of Proposition \ref{ref-5.1-45}]
The forgetful functor
\[
\Dif(\AAA,C(\ccc))
\lra C(\ccc): F
\longmapsto F(A)
\]
  has a left adjoint
\[
L_A: C(\ccc) \lra \Dif(\AAA,C(\ccc)):
C
\longmapsto (A'
\longmapsto \AAA(A,A') \otimes_k C).
\]
Since every $\AAA(A,A')$ has projective components, $L_A$ preserves
monomorphisms. Since  $\AAA(A,A')\otimes-$ preserves acyclic complexes
$L_A$
preserves weak equivalences and hence $L_A$ preserves trivial
cofibrations. It now easily follows that if $F$ has the lifting
property with respect to trivial cofibrations then so does every $F(A)$.

Hence (2)
follows if we prove (1).
For (1), we use Proposition \ref{ref-5.2-46}  for
$$H^i: \Dif(\AAA,C(\ccc)) \lra \ccc^{|\AAA|}: M \longmapsto
(H^iM(A))_{A}.$$
To see that $(H^i)_i$ is effaceable, we can take for $A \in
\Dif(\AAA,\ccc)$ the monomorphism $A \lra \mathrm{cone}(1_A)$.
\end{proof}
\subsection{The derived Gabriel-Popescu theorem}
\label{ref-5.2-47}
The results in this section are probably well-known but we have not
been able to locate a reference.

Let $j: \UUU \lra \ccc$
be a $k$-linear functor from a small $k$-linear category inducing a
\emph{localization}
\begin{equation}
\label{ref-5.1-48}
c: \ccc \lra \Mod(\UUU): C \longmapsto
\ccc(j(-),C)
\end{equation}
  By this we mean that $c$ is fully
faithful and has an exact left adjoint. By the Gabriel-Popescu
theorem \cite{GP} $c$ is a localization if $j$ is fully faithful and
generating. However in our applications to ringed spaces
$j$ will \emph{not} be fully faithful.  Necessary and sufficient
conditions for $c$ to be a localization were given in \cite{lowen1}.

We start with the following easy
result.
\begin{theorems}
\label{ref-5.2.1-49} The functor
\[
D(\Cscr)\r D(\UUU)
\]
which sends a fibrant object $A\in C(\Cscr)$ to
$\Hom_\Cscr(j(-),A)$
preserves $\RHom$.
\end{theorems}
\begin{proof}
We need to prove that for fibrant $A$, $B$ we have a quasi-isomorphism
\[
\Hom_\Cscr(A,B)=\RHom_\UUU(\Hom_\Cscr(j(-),A),\Hom_\Cscr(j(-),B))
\]
Since $c$ is fully faithful we have
\[
\Hom_\Cscr(A,B)=\Hom_\UUU(\Hom_\Cscr(j(-),A),\Hom_\Cscr(j(-),B))
\]
And since $c$ has an exact left adjoint we easily deduce that
$\Hom_\UUU(-,\Hom_\Cscr(j(-),B))$ preserves acyclic complexes. Thus
\[
\Hom_\UUU(\Hom_\Cscr(j(-),A),\Hom_\Cscr(j(-),B))=\RHom_\UUU(\Hom_\Cscr(j(-),A),\Hom_\Cscr(j(-),B))
\]
This finishes the proof.
\end{proof}
Now we discuss a more sophisticated derived version of the
Gabriel-Popescu theorem.
Let $l:\FFF\r C(\ccc)$ be any fully faithful
DG-functor such that:
\begin{enumerate}
\item $l(\FFF)$ consists of fibrant complexes;
\item every object in $j(\UUU)$ is quasi-isomorphic to an object in $l(\FFF)$;
\item the only cohomology of an object in $l(\FFF)$ is in degree zero 
and lies in
$j(\UUU)$.
\end{enumerate}
For example $\FFF$ could consist of  injective
resolutions for the objects $j(U)$.
\begin{theorems} \label{ref-5.2.2-50} The functor
\[
D(\Cscr)\r D(\FFF)
\]
which sends a fibrant object $A\in C(\Cscr)$ to
$\Hom_\Cscr(l(-),A)$
preserves $\RHom$.
\end{theorems}
The rest of this section will be devoted to the proof of this
theorem. Along the way we introduce some notations which will also be used
afterwards.

To start we note the following.
\begin{lemmas} It is sufficient to prove Theorem \ref{ref-5.2.2-50} for one
   particular choice of~$\frak{\FFF}$.
\end{lemmas}
\begin{proof}
Let $\FFF_1$ be the full DG-subcategory of $C(\ccc)$ of all fibrant complexes
satisfying (3). Then $\FFF \lra \FFF_1$ is a
quasi-equivalence, hence the corresponding functor $\Dif(\FFF_1)\r
   \Dif(\FFF)$ preserves $\RHom$ (Lemma \ref{ref-2.2.2-7}). Therefore, if
   Theorem \ref{ref-5.2.2-50} is true for one $\frak{\FFF}$, it is true for
   $\frak{\FFF}_1$ and then it is true for all~$\frak{\FFF}$.
\end{proof}
Ideally we would want to choose $\frak{f}$ in such a way that there is
a corresponding DG-functor $\frak{\UUU}\r \frak{\FFF}$.  It is not
clear to us that this is possible in general. However, as we now 
show, it is possible after
replacing $\UUU$ by a $k$-cofibrant resolution.

Let $r: \ovl{\UUU} \lra \UUU$ be a $k$-cofibrant resolution of $\UUU$.
Let $i: \ccc \lra C(\ccc)$ be the canonical inclusion and let $ijr
\lra E$ be a fibrant resolution for the model-structure on
$\Dif(\ovl{\UUU}^{\op},\ccc)$ of \S\ref{ref-5.1-44}. By Proposition
\ref{ref-5.1-45} this yields fibrant replacements $ij(U) \lra E(U)$
natural in $U \in \ovl{\UUU}$. We define a new DG-category
$\ovl{\ovl{\UUU}}$ with the same objects as $\UUU$ and
\[
\ovl{\ovl{\UUU}}(U,V) = \Hom_{\ccc}(E(U),E(V))
\]
By construction we have the following commutative diagram of DG-functors
\[
\begin{CD}
\UUU @>ij >> C(\Cscr)\\
@A r AA         @AA E(-) A\\
\bar{\UUU} @>>f> \ovl{\ovl{\UUU}}
\end{CD}
\]
where $r$ is a quasi-equivalence and $E(-)$ is fully faithful.
\begin{proof}[Proof of Theorem \ref{ref-5.2.2-50}]
We will prove the theorem with $\FFF=\ovl{\ovl{\UUU}}$.

\begin{step} Let $A\in C(\Cscr)$ be fibrant.
We first claim that the canonical map
\[
\Hom_\Cscr(E(-),A)\r Rf^!f_\ast \Hom_\Cscr(E(-),A)
\]
is a quasi-isomorphism in $\Dif(\ovl{\ovl{\UUU}})$.
The proof is based on the following computation for $U\in \Ob(\UUU)$
\begin{align*}
\Hom_\Cscr(E(U),A)
&\cong \RHom_\UUU(\Hom_\Cscr(j(-),E(U)),\Hom_\Cscr(j(-),A))\\
&\cong\RHom_{\ovl{\UUU}}(\Hom_\Cscr(jr(-),E(U)),\Hom_\Cscr(jr(-),A))\\
&\cong\RHom_{\ovl{\UUU}}(\ovl{\ovl{\UUU}}(-,U),\Hom_\Cscr(E(-),A))\\
&= (Rf^!f_\ast \Hom_\Cscr(E(-),A))(U)
\end{align*}
The first line is Theorem \ref{ref-5.2.1-49}. In the second
line we use the fact that $\ovl{\UUU}\r \UUU$ is a quasi-equivalence
together with Lemma \ref{ref-2.2.2-7}. The third line is a change of notation
and the fourth line is an easy verification.
\end{step}
\begin{step} Now we finish the proof of the theorem. We compute for fibrant
$A,B\in C(\Cscr)$
\begin{align*}
\RHom_{\ovl{\ovl{\UUU}}}(\Hom_\Cscr(E(-),A)&,\Hom_\Cscr(E(-),B))\\
&\cong
\RHom_{\ovl{\ovl{\UUU}}}(\Hom_\Cscr(E(-),A),Rf^!f_\ast\Hom_\Cscr(E(-),B))\\
&\cong\RHom_{\ovl{\UUU}}(f_\ast \Hom_\Cscr(E(-),A),f_\ast\Hom_\Cscr(E(-),B))\\
&\cong\RHom_{\ovl{\UUU}}(\Hom_\Cscr(j(-),A),\Hom_\Cscr(j(-),B))\\
&\cong\RHom_{{\UUU}}(\Hom_\Cscr(j(-),A),\Hom_\Cscr(j(-),B))\\
&=\Hom_\Cscr(A,B)
\end{align*}
where we have once again used Theorem \ref{ref-5.2.1-49} and
the fact that $\ovl{\UUU}\r \UUU$ is a quasi-equivalence. \qed
\end{step}
\def\qed{}\end{proof}
\begin{remarks}
\label{ref-5.2.4-51} If $\frak{u}$ is already $k$-cofibrant then we may take
$\bar{\UUU}=\UUU$. By letting $E$ be an injective resolution of $ij$
we obtain injective resolutions $U\r E(U)$ of $j(U)$ natural in $U$.
\end{remarks}
\subsection{Hochschild complexes}
\label{ref-5.3-52}
Let $\III = \Inj(\ccc)$ be the category of injectives in $\ccc$ and
let $l:\FFF\r C(\Cscr)$ be as in \S\ref{ref-5.2-47}. In this section 
we prove the following
comparison result.
\begin{theorems}\label{ref-5.3.1-53}
There is a quasi-isomorphism $\CC_{\mathrm{sh}}(\III) \cong
\CC_{\mathrm{sh}}(\FFF)$.
\end{theorems}
\begin{corollarys}
\label{ref-5.3.2-54}
  $\CC_{\mathrm{sh}}(\III)$ has small cohomology.
\end{corollarys}
\begin{proof} It is clear that we may take $\FFF$ to be small.
\end{proof}
\begin{proof}[Proof of Theorem \ref{ref-5.3.1-53}]
We define the ${\III}{-}\FFF$-bimodule
\[
X(U,E) =\Hom_{\ccc}(l(U),E).
\]
By Lemma \ref{ref-2.2.2-7} and the derived Gabriel-Popescu Theorem 
(Theorem \ref{ref-5.2.2-50})
\begin{align*}
\RHom_{\FFF}({X(-,E)},{X(-,F)})
&\cong \Hom_{\ccc}(E,F)\\
&\cong {\III}(E,F)
\end{align*}
On the other hand,
\begin{align*}
\RHom_{{\III}^{\op}}(X(V,-),X(U,-)) &\cong
\RHom_{{\III}^{\op}}(\ccc(l(V),-),\ccc(l(U),-))\\
&\cong \RHom_{\ccc}(l(U),l(V))\\
&\cong \FFF(U,V)
\end{align*}
where the first line is a consequence of Lemma \ref{ref-5.3.3-55} below.
The result now follows from 
Proposition \ref{ref-4.3.3-36}.
\end{proof}
\begin{lemmas}\label{ref-5.3.3-55}
Suppose a small abelian category $\aaa$ has enough injectives in $\add(\JJJ)$
for $\JJJ
\subset \Inj(\aaa)$ (i.e., for every object in $\aaa$ there is a mono into a
finite sum of injectives in $\JJJ$). Consider
$\aaa^{\op}
\lra
\Mod(\JJJ): A
\longmapsto
\aaa(A,-)$. For
$A,B \in \aaa$, we have $$\RHom_{\aaa}(A,B) \cong \RHom_{\JJJ^{\op}}(\aaa(B,-),
\aaa(A,-)).$$
\begin{proof}
Let $B \lra I^{\cdot}$ be an injective resolution of $B$ in $\add(\JJJ)$.
Then $\aaa(I^{\cdot},-)\lra \aaa(B,-)$ is a resolution in $\Mod(\JJJ)$,
and every object
$\aaa(I^i,-) = \aaa(\oplus_{k=1}^nJ^i_k,-) = \oplus_{k = 1}^n\JJJ(J^i_k,-)$ is
projective.
\end{proof}
\end{lemmas}
In the proof of Theorem \ref{ref-5.3.1-53} we have used this lemma in the
case
$\JJJ=\Inj(\Ascr)$. The added generality will be used in \S
\ref{ref-7.7-97}.

\subsection{A spectral sequence}
We keep the same notations as in \S\ref{ref-5.2-47}. The main theorem of
this section is an  interesting spectral sequence which
relates the Hochschild cohomology of $\UUU$ to that of
$\III=\Inj(\Cscr)$.
\begin{theorems} \label{ref-5.4.1-56}
There is a convergent, first quadrant spectral sequence
\begin{equation}
\label{ref-5.2-57}
E_2^{pq}:\mathrm{H}\mathbf{C}_{\mathrm{sh}}^p(\UUU,\Ext^q_{\Cscr}(j(-),j(-))
\Rightarrow \mathrm{H}\mathbf{C}_{\mathrm{sh}}^{p+q}(\III)
\end{equation}
\end{theorems}
The proof of Theorem \ref{ref-5.4.1-56} depends on the following
technical result.
\begin{lemmas} \label{ref-5.4.2-58} There is a (non 
$B_\infty$-)quasi-isomorphism
$\CC(\ovl{\UUU},\ovl{\ovl{\UUU}})
\cong \CC_{\mathrm{sh}}(\frak{\III})$.
\end{lemmas}
\begin{proof}
By Theorem \ref{ref-5.3.1-53} it is sufficient to construct a
quasi-isomorphism $\CC(\ovl{\UUU},\ovl{\ovl{\UUU}})\cong
\CC_{\mathrm{sh}}(\ovl{\ovl{\UUU}})$.
   Let $\VVV$ be a cofibrant
resolution of $\ovl{\ovl{\UUU}}$ as a
$\ovl{\ovl{\UUU}}{-}\ovl{\UUU}$-bimodule.

Of course
\begin{equation}
\label{ref-5.3-59}
\RHom_{\ovl{\ovl{\UUU}}^{\op}}(\VVV(V,-),\VVV(U,-))
\cong 
\RHom_{\ovl{\ovl{\UUU}}^{\op}}(\ovl{\ovl{\UUU}}(V,-),\ovl{\ovl{\UUU}}(U,-))
\cong \ovl{\ovl{\UUU}}(U,V)
\end{equation}
  so we compute
\begin{align*}
\RHom_{\ovl{\UUU}}(\VVV(-,U),\VVV(-,V))
&\cong \RHom_{\ovl{\UUU}}(\ovl{\ovl{\UUU}}(-,U),\ovl{\ovl{\UUU}}(-,V))\\
&\cong \RHom_{\ovl{\UUU}}(\Hom_{\ccc}(jr(-),E(U)),\Hom_{\ccc}(jr(-),E(V)))\\
&\cong \RHom_{\UUU}(\Hom_{\ccc}(j(-),E(U)),\Hom_{\ccc}(j(-),E(V)))\\
&\cong \Hom_{\ccc}(E(U),E(V))\\
&\cong \ovl{\ovl{\UUU}}(U,V)
\end{align*}
where we have use that $\ovl{\UUU}\r \UUU$ is a quasi-isomorphism,
together with Theorem \ref{ref-5.2.1-49}.
  The result now follows from Theorem \ref{ref-4.1.1-25}.1.
\end{proof}

\begin{proof}[Proof or Theorem \ref{ref-5.4.1-56}]  We use Lemma
   \ref{ref-5.4.2-58}. Let $\bar{\UUU}$ be a semi-free \cite[\S13.4]{Drinfeld2}
   resolution of~$\UUU$. This is in particular a $k$-cofibrant
   resolution $\bar{\UUU}\r \UUU$ concentrated in non-positive degree.
   The latter implies that the truncations $\tau_{\le
     n}\ovl{\ovl{\UUU}}$ and $\tau_{\ge
     n}\ovl{\ovl{\UUU}}=\ovl{\ovl{\UUU}}/\tau_{< n}\ovl{\ovl{\UUU}}$
   are $\bar{\UUU}$-bimodules. Recall that the homology of
   $\ovl{\ovl{\UUU}}$ is $\Ext^*_\Cscr(j(-),j(-))$ so it lives in
   non-negative degree. So up to quasi-isomorphism we may replace
   $\ovl{\ovl{\UUU}}$ by $\WWW=\tau_{\ge 0}\ovl{\ovl{\UUU}}$.  We put
   the ascending filtration $(\tau_{\le n}\WWW)_n$ on $\WWW$.  This
   filtration is positive since $\tau_{<0}\WWW=0$.

We claim that the obvious map
\[
\bigcup_n \mathbf{C}(\bar{\UUU},\tau_{\le n} \WWW)
\lra\mathbf{C}(\bar{\UUU},\WWW)
\]
is a quasi-isomorphism. To prove this it is sufficient to show that
for a fixed $i$ the map
\[
  \mathrm{H}\mathbf{C}^i(\bar{\UUU},\tau_{\le n} \WWW)
\lra\mathrm{H}\mathbf{C}^i(\bar{\UUU},\WWW)
\]
is a quasi-isomorphism for large $n$. Equivalently by the long exact
sequence for Hochschild cohomology,
$\mathrm{H}\mathbf{C}^i(\bar{\UUU},\tau_{\ge n} \WWW)$ should
be zero for large $n$. Now by \eqref{ref-2.5-15} we have
\[
\mathrm{H}\mathbf{C}^i(\bar{\UUU},\tau_{\ge n} \WWW)
=\RHom_{\bar{\UUU}^\circ\otimes \bar{\UUU}}(\bar{\UUU},\tau_{\ge n} \WWW[i])
\]
Since a cofibrant $\bar{\UUU}^{\op}\otimes \bar{u}$-resolution of $\bar{\UUU}$
may also be chosen to live in non-positive degree it is clear that this
is zero for $n>i$.

So the spectral sequence associated to the filtered complex
$\bigcup_n \mathbf{C}(\bar{\UUU},\tau_{\le n} \WWW)$
converges to $\mathrm{H}\mathbf{C}^\ast_\ast(\III)$. The associated
graded complex is
\[
\bigoplus_n \mathbf{C}(\bar{\UUU},\Ext_\Cscr^n(j(-),j(-))[-n])
\]
and the homology of this graded complex is:
\[
\bigoplus_{mn}\mathrm{H}\mathbf{C}_{\mathrm{sh}}^{m-n}(\UUU,\Ext_\Cscr^n(j(-),j(-)))
\]
After the appropriate reindexing we obtain the desired result.
\end{proof}
\subsection{Application of a censoring relation}
Lemma \ref{ref-5.4.2-58} has the following useful variant:
\begin{propositions}
\label{ref-5.5.1-60}
Assume that $\UUU$ is equipped with a (possibly trivial)
censoring relation $\Rscr$ (see \S\ref{ref-4.3-32}) such that
\[
\Ext^i_\Cscr(j(U),j(V))=0\text{ for }i>0\text{ and }(U,V)\in \Rscr
\]
Then there is a quasi-isomorphism
\[
\CC_{\mathrm{ab}}(\Cscr)\cong \CC_{\mathrm{sh}}(\UUU)
\]
\end{propositions}
\begin{proof} We define the ${\III}{-}\UUU$-bimodule
\[
X(U,E) =\Hom_{\ccc}(j(U),E).
\]
By Theorem \ref{ref-5.2.1-49}
\begin{align*}
\RHom_{\UUU}({X(-,E)},{X(-,F)})
&\cong \Hom_{\ccc}(E,F)\\
&\cong {\III}(E,F)
\end{align*}
On the other hand for $(U,V)\in \Rscr$
\begin{align*}
\RHom_{{\III}^{\op}}(X(V,-),X(U,-)) &\cong
\RHom_{{\III}^{\op}}(\ccc(j(V),-),\ccc(j(U),-))\\
&\cong \RHom_{\ccc}(j(U),j(V))\\
&\cong \UUU(U,V)
\end{align*}
where the first line is a consequence of Lemma \ref{ref-5.3.3-55}.
The result now follows from Proposition \ref{ref-4.3.3-36}.
\end{proof}

\section{Basic results about Hochschild cohomology of
   abelian categories}
\label{ref-6-61}
Let $\Ascr$ be a small abelian category. By definition we have
$\CC_{\mathrm{ab}}(\Ascr)=\CC_{\mathrm{sh}}(\III)$ with $\III=\Inj\,
\Ind (\AAA)$ (see
\S\ref{ref-3-18}).  The embedding $\aaa \lra \Ind(\aaa)$ satisfies the
hypotheses on $j$ in \S\ref{ref-5-43} so the results of that section
apply. We will now translate them to the current setting.

The first result below relates the Hochschild cohomology of $\Ascr$ to
that of suitable small DG-categories.  Let ${}^e\!\!\Ascr$ be the full
DG-subcategory of $C(\Ind(\aaa))$ spanned by all positively graded
complexes of injectives whose only cohomology is in degree zero and
lies in $\aaa$ and let ${}^e\! D^b(\Ascr)$ be spanned by all left
bounded complexes of injectives with bounded cohomology in~$\aaa$.
Note that for example by \cite[Prop.\ 2.14]{lowenvdb1} we have
\[
H^0({}^e \!D^b(\Ascr))\cong D^b(\Ascr)
\]
so ${}^e \!D^b(\Ascr)$ is a DG-enhancent for $D^b(\Ascr)$.
\begin{theorem}\label{ref-6.1-62} There are quasi-isomorphisms
\[
\CC_{\mathrm{ab}}(\aaa)
\cong \CC_{\mathrm{sh}}({}^e\!\!\Ascr) \cong
\CC_{\mathrm{sh}}({}^e\!D^b(\Ascr))
\]
\begin{proof}
The first quasi-isomorphism follows
from Theorem \ref{ref-5.3.1-53}. For the second quasi-isomorphism put
$\AAA={}^e\Ascr$ and let $\BBB$ be the
closure of $\AAA$ in $C(\Ind\Ascr))$ under finite cones and
shifts. Then $\BBB\r {}^e D^b(\Ascr)$ is a
quasi-equivalence. Furthermore the functor $B \longmapsto
\Hom_{C(\Ind\Ascr)}(-,B)$ defines an embedding $\BBB\r \Dif(\AAA)$
whose image consists of cofibrant objects. We may now invoke Theorem 
\ref{ref-4.4.1-39}.
\end{proof}
\end{theorem}
In \cite{keller6}, Bernhard Keller defines the Hochschild complex
$\CC_{\mathrm{ex}}(\eee)$ of an exact category $\eee$ as
$\CC_{\mathrm{sh}}(\qqq)$ for a DG-quotient $\qqq$ of $Ac^b(\eee) \lra
C^b(\eee)$, where $C^b(\eee)$ is the DG-category of bounded complexes
of $\eee$-objects and $Ac^b(\eee)$ is its full DG-subcategory of
acyclic complexes. We endow the abelian category $\aaa$ with the exact
structure given by all exact sequences.
\begin{theorem}\label{ref-6.2-63} There is a quasi-isomorphism
$$\CC_{\mathrm{ab}}(\aaa) \cong \CC_{\mathrm{ex}}(\aaa).$$
\begin{proof}
This follows from Theorem \ref{ref-6.1-62} and  Lemma \ref{ref-6.3-64}
below.
\end{proof}
\end{theorem}
\begin{lemma}\label{ref-6.3-64}
${}^e D^b(\Ascr)$ is a DG-quotient of
$Ac^b(\aaa) \lra C^b(\aaa)$.
\begin{proof}
We sketch the proof.
Let $\CCC$ be the
following DG-category: the objects of
$\CCC$ are quasi-isomorphisms $f: C \lra I$ with $C \in C^b(\aaa)$ and $I \in
{}^e D^b(\Ascr)$. Morphisms from $f$ to $g:
D
\lra J$ are maps
$\mathrm{cone}(f)
\lra
\mathrm{cone}(g)$ with zero component $I \lra D[1]$. The two
projections yield a diagram $C^b(\aaa)
\longleftarrow \CCC \lra {}^e D^b(\Ascr)$, for which it
is easily seen that
$C^b(\aaa) \longleftarrow \CCC$ is a quasi-equivalence
and $\CCC \lra {}^e D^b(\Ascr)$ induces the
exact sequence of associated triangulated categories $H^0(Ac^b(\aaa)) \lra
H^0(C^b(\aaa)) \lra D^b_{\aaa}(\Ind(\aaa)) \cong D^b(\aaa)$.
By \cite{kelcyex,Drinfeld2}, this proves the statement.
\end{proof}
\end{lemma}
\begin{remark} Since $\CC_{\mathrm{ex}}$ is easily seen to satisfy
   $\CC_{\mathrm{ex}}(\eee) \cong \CC_{\mathrm{ex}}(\eee^{\op})$ (in $D(k)$) if
   $\eee^{\op}$ is endowed with the opposite exact sequences of $\eee$,
   by Theorem \ref{ref-6.2-63} we have in particular that
   $\CC_{\mathrm{ab}}(\aaa) \cong \CC_{\mathrm{ab}}(\aaa^{\op})$. It is
   a pleasant excercise to derive this result directly from our
   definition.
\end{remark}

The following result, of theoretical interest, is a restatement of a
special case of the
spectral sequence \eqref{ref-5.2-57}. It compares the
Hochschild cohomology of $\Ascr$ as linear and as abelian category.
\begin{proposition} \label{ref-6.5-65}
There is a convergent, first quadrant spectral sequence
\begin{equation}
\label{ref-6.1-66}
E_2^{pq}:\mathrm{H}\mathbf{C}_{\mathrm{sh}}^p(\Ascr,\Ext^q_{\Ascr}(-,-))
\Rightarrow \mathrm{H}\mathbf{C}_{\mathrm{ab}}^{p+q}(\Ascr)
\end{equation}
\end{proposition}
\begin{proof} We only need to remark that Yoneda-$\Ext$ computed in
   $\Ascr$ and $\Ind \Ascr$ is the same. See for example \cite[Prop.\
   2.14]{lowenvdb1}.
\end{proof}
The following result shows that if $\Ascr$ has enough injectives then
there is no need to pass to $\Ind \Ascr$.

\begin{theorem}\label{ref-6.6-67} Assume that $\Ascr$ has enough
   injectives and put $\III=\Inj \Ascr$.
There is a quasi-isomorphism
\[
\CC_{\mathrm{ab}}(\aaa)  \cong \CC_{\mathrm{sh}}(\III).
\]
\end{theorem}
\begin{proof}
   Let $\KKK$ be the the full DG-subcategory of $C(\Ascr)$ spanned by
   all positively graded complexes of $\III$-objects whose only
   cohomology is in degree zero. Then the inclusion $\KKK\r
   {}^e\!\! \Ascr$ is a quasi-equivalence
   so using Theorem \ref{ref-6.1-62} it is sufficient to show that $\KKK$ and
   $\III$ have isomorphic Hochschild complexes. We embed $\KKK$ in
   $\Dif(\III)$ via the functor $E\mapsto \Hom_\Cscr(-,E)$. Then $\KKK$
   is mapped to right bounded projective complexes. Such complexes are
   cofibrant and hence we may use Theorem \ref{ref-4.4.1-39} to deduce
   $\CC_{\mathrm{sh}}(\KKK)\cong \CC_{\mathrm{sh}}(\III$).
\end{proof}
The following corollary will be used in \S \ref{ref-7.7-97}.
\begin{corollary}
\label{ref-6.7-68}
Assume that $\Ascr$ has enough
   injectives in $\add(\JJJ)$. There is a quasi-isomorphism
\[
\CC_{\mathrm{ab}}(\aaa)  \cong \CC_{\mathrm{sh}}(\JJJ).
\]
\end{corollary}
\begin{proof}
Put $\III = \Inj(\aaa)$. By Lemma \ref{ref-5.3.3-55}, the inclusion $\JJJ \lra
\III$ satisfies the hypotheses of Proposition \ref{ref-4.3.4-38}, hence
$\CC_{\mathrm{sh}}(\III)  \cong \CC_{\mathrm{sh}}(\JJJ)$. The result now
follows from Theorem \ref{ref-6.6-67}.
\end{proof}

\begin{corollary}
\label{ref-6.8-69}
If $\aaa$ is a  small abelian category then there is a quasi-isomorphism
$$\CC_{\mathrm{ab}}(\aaa) \cong \CC_{\mathrm{ab}}(\Ind(\aaa)).$$
\end{corollary}
\begin{proof} Immediate from Theorem \ref{ref-6.6-67} and the definition.
\end{proof}
Theorem \ref{ref-6.6-67} applies in particular if $\Ascr$ is a
Grothendieck category. So the results in \S\ref{ref-5-43} (with
$\Cscr$ replaced by $\Ascr$) may be reinterpreted as being about the
Hochschild complex of a Grothendieck category. We mention in
particular Theorem \ref{ref-5.3.1-53} and Proposition \ref{ref-5.5.1-60} which
shows how to compute $\CC_{\mathrm{ab}}(\aaa)$ in terms of generators,
Corollary \ref{ref-5.3.2-54} which shows that $\CC_{\mathrm{ab}}(\aaa)$ has
small homology and the spectral sequence \eqref{ref-5.2-57}
abutting to $\mathrm{H}\CC^\ast_{\mathrm{ab}}(\aaa)$.

The following corollary to
Proposition \ref{ref-5.5.1-60} was our original motivation
for starting this project.
\begin{corollary}
\label{ref-6.9-70}
Let $\AAA$ be a small $k$-category. There is a quasi-isomorphism
\[
\CC_{\mathrm{ab}}(\Mod(\AAA)) \cong \CC_{\mathrm{sh}}(\AAA).
\]
In particular, for a $k$-algebra $A$, there is a quasi-isomorphism
\[
\CC_{\mathrm{ab}}(\Mod(A)) \cong \CC_{\mathrm{sh}}(A).
\]
\end{corollary}
\begin{proof}
We apply Proposition \ref{ref-5.5.1-60} with  the Yoneda embedding
$j:\AAA\r \Mod(\AAA)$.
\end{proof}
The following proposition shows that $H\CC^\ast_{\mathrm{ab}}(\Ascr)$
defines elements in the center of $D^b(\Ascr)$.
\begin{proposition}
\label{ref-6.10-71}
There is a homomorphism of graded rings
\[
\sigma_\Ascr:H\CC_{\mathrm{ab}}^\ast(\Ascr)\r Z(D^b(\Ascr))
\]
where on the right hand side we view $D^b(\Ascr)$ as a \emph{graded}
category in the usual way.
\end{proposition}
\begin{proof}  Put
   $\AAA={}^e \!D^b(\Ascr)$. Then by
   \S\ref{ref-4.5-41} there is a homomorphism of graded rings
\[
\sigma_\Ascr:H\CC^\ast_{\mathrm{sh}}(\AAA)\r Z(H^\ast(\AAA))
\]
We define $\sigma_\Ascr$ as the composition
\[
H\CC^\ast_{\mathrm{ab}}(\aaa)\cong H\CC^\ast_{\mathrm{sh}}(\AAA)
\xrightarrow{\sigma_\AAA}Z(H^\ast(\AAA))\cong Z(D^b(\AAA))
\]
where the first isomorphism comes from Theorem \ref{ref-6.1-62}.
\end{proof}
\begin{remark} We may think of $\sigma_\Ascr$ as defining
   ``universal'' elements in $\Ext^\ast_\Ascr(M,M)$ for every $M\in
   D^b(\Ascr)$. These universal elements are closely related to Atiyah
   classes in algebraic geometry.  See for example~\cite{BuchweitzFlenner}.
\end{remark}
\section{Hochschild cohomology for ringed spaces and schemes}
\label{ref-7-72}
\subsection{Discussion and statement of the main results}
Below let $(X,\Oscr)$ be a $k$-linear possibly non-commutative ringed space.
We define the
Hochschild complex of $X$ as
\[
\CC(X)=\CC_{\mathrm{ab}}(\Mod(X))
\]
where $\Mod(X)$ is the category of sheaves of right modules over
$X$.  For the purpose of clarity we will sometimes use the notation
  $\CC(X,\Oscr)$  for $\CC(X)$. Note
that in the definition of $\CC(X)$ the bimodule structure of
$\Oscr$ does not enter explicitly.

As $\Mod(X)$ has enough injectives an equivalent definition
(using Theorem \ref{ref-6.6-67}) for $\CC(X)$ would be
\[
\CC(X)=\CC_{\mathrm{sh}}(\Inj\Mod(X))
\]
Recall that $\CC(X)$ describes the deformation theory of the abelian
category $\Mod(X)$ (as explained in \S\ref{ref-3-18}) but \emph{not} of the
ringed space $(X,\Oscr)$. This is a related but different
problem.\footnote{To obtain the correct correspondence we should
   deform $(X,\Oscr)$ in the category of algebroids over $X$, see
   \cite{Ko4}}

We now summarize some of the results we will prove about $\CC(X)$.  We
would like to think of $H\CC^\ast(X)$ as defining a (generalized)
cohomology theory for ringed spaces. A first indication for this is
that $\CC(-)$ is a contravariant functor on open embeddings of ringed
spaces and associated to an open covering $X=U\cup V$ there is a
corresponding Mayer-Vietoris long exact sequence (see \S\ref{ref-7.9-103})
\begin{equation}
\label{ref-7.1-73}
\cdots \rightarrow H\CC^{i-1}(U\cap V)\rightarrow
  H\CC^i(X)\rightarrow H\CC^i(U)\oplus H\CC^i(V)
\rightarrow H\CC^i(U\cap V)\rightarrow \cdots
\end{equation}
Let $\underline{k}$ be the constant sheaf with values in $k$. Our
next interesting result is an isomorphism between Hochschild cohomology
and ordinary cohomology
\begin{equation}
\label{ref-7.2-74}
H\CC^\ast(X,\underline{k}) \cong H^\ast(X,\underline{k}).
\end{equation}
  In \cite{Baues} Baues
shows that the singular cochain complex of a topological space is a
$B_\infty$-algebra. Thus \eqref{ref-7.2-74} suggests that
$\CC(X,\underline{\ZZ})$ should be viewed as an algebraic analog of
the singular cochain complex of $X$.

\medskip

Now we discuss some more specific results for Hochschild cohomology.
For any subposet $\uuu$ of $\Open(X)$ let $\UUU=\UUU({\uuu})$ be the
linear category with $\Ob(\UUU) = \uuu$ and
\[
\UUU(U,V) = \begin{cases} \ooo(U)& \text{if}\,\, U \subset V\\
0&\text{otherwise}
\end{cases}\]
First assume that $\bbb$ is a basis of $X$ of \emph{acyclic opens},
i.e.\ for $U \in \bbb$: $H^i(U,\ooo_U) = 0$ for $i>0$.
Put $\BBB = \UUU({\bbb})$. Our first result (see \S\ref{ref-7.3-82})
is that there is a quasi-isomorphism
\begin{equation}
\label{ref-7.3-75}
\CC(X) \cong \CC_{\mathrm{sh}}({\BBB}).
\end{equation}
In \cite{GS1,GS2} Gerstenhaber and Schack define the 
$k$-\emph{relative} Hochschild
complex of a presheaf of rings. In order to make a connection with 
our setting let us
assume that $k$ is a field. Let $\Bscr$ be an acyclic basis as above
and let $\Oscr_\Bscr$ be the restriction of $\Oscr$ to $\Bscr$,
considered as a presheaf of rings. It is implied in \cite{GS1,GS2}
that
\begin{equation}
\label{ref-7.4-76}
\CC_{GS}(\Oscr_\Bscr)\overset{\mathrm{def}}{=}\RHom_{\Oscr_\Bscr^{\op}\otimes
\Oscr_\Bscr}(\Oscr_\Bscr,\Oscr_\Bscr)
\end{equation}
is a reasonable definition for the Hochschild complex of $X$. We will
show that this is
true. Indeed it follows from combining \eqref{ref-7.3-75} with
\eqref{ref-7.8-80}\eqref{ref-7.9-81} below that
\[
\CC(X)\cong \CC_{GS}(\Oscr_\Bscr)
\]
Let $k$ be general again. We now specialize to the case where
$X$ is a quasi-compact separated scheme over~$k$.

  Let $X =
\cup_{i=1}^nA_i$  be  a finite affine open covering of $X$ and
  let $\Ascr$ be
the closure of this covering under intersections. Put
$\AAA=\UUU(\Ascr)$. Of course $\AAA$ is not a basis for $X$ but
nevertheless
we have the following analog of
\eqref{ref-7.3-75} (see \S\ref{ref-7.5-87})
\begin{equation}
\label{ref-7.5-77}
\CC(X) \cong \CC_{\mathrm{sh}}({\AAA}).
\end{equation}
Let $\Qch(X)$ be the category of quasi-coherent
$\Oscr$-modules. We will prove (see \S\ref{ref-7.7-97})
\[
\CC(X) \cong \CC_{\mathrm{ab}}(\Qch(X)).
\]
and if $X$ is noetherian we even have
\[
\CC(X) \cong \CC_{\mathrm{ab}}(\coh(X)).
\]
where  $\coh(X)$
is the category of coherent $\Oscr$-modules.

If $k$ is a field then in \cite{swan}, Richard G. Swan defines the
Hochschild complex of $(X,\Oscr)$ to be
\[
\CC_{\mathrm{Swan}}(X)\overset{\text{def}}{=} \RHom_{X\times
   X}(\ooo_\Delta,\ooo_\Delta)
\]
where $\Delta\subset X\times X$ is the diagonal. We prove
\begin{equation}
\label{ref-7.6-78}
\CC(X)\cong \CC_{\mathrm{Swan}} (X)
\end{equation}
This is known in the finite type case since in that case
by
\cite[\S 3]{swan}
\begin{equation}
\label{ref-7.7-79}
\CC_{\mathrm{Swan}}(X)\cong \CC_{\mathrm{GS}} (\Oscr_\Bscr)
\end{equation}
where $\Bscr$ is  the acyclic basis of all afine opens.
\subsection{Presheaves of modules over presheaves of rings}
Let $\Bscr$ be a poset and let
$\Oscr$ be a presheaf of
$k$-algebras on $\Bscr$.  Let $\Pre(\Oscr)$ be the category of
$\Oscr$-modules.

Associated to $\Oscr$ is a small $k$-linear category $\BBB$ with
$\Ob(\BBB) =
\Bscr$ and
\[
\BBB(U,V) = \begin{cases} \Oscr(U) &\text{if }\,U \subset V \\ 0 & 
\text{otherwise}
\end{cases}
\]
\begin{propositions}
There is a
quasi-isomorphism
\begin{equation}
\label{ref-7.8-80}
\CC_{\mathrm{ab}}(\Pre(\Oscr)) \cong \CC_{\mathrm{sh}}(\BBB).
\end{equation}
\end{propositions}
\begin{proof}
For $V\in \Bscr$ let $P_V=\BBB(-,V)$ be the extension by zero of
$\Oscr\mid V$.
The $(P_V)_V$ form a system of small projective generators for
$\Pre(\Oscr)$ and $\BBB \lra \Pre(\Oscr): V \longmapsto P_V$ yields an
equivalence of categories
\[
\Mod(\BBB) \cong \Pre(\Oscr).
\]
Hence the result is just a rephrazing of Corollary \ref{ref-6.9-70}.
\end{proof}
We now discuss the relation with the papers \cite{GS1,GS2} by
Gerstenhaber-Schack. These authors work with relative Hochschild
cohomology which makes it somewhat difficult to translate their
results to our situation.
So for simplicity we assume that $k$ is a
field.  The Hochschild complex of $\Oscr$ according to
\cite{GS1,GS2} is
\[
\CC_{\mathrm{GS}}(\Oscr)\overset{\text{def}}{=}\RHom_{\Oscr^{\op}\otimes
   \Oscr}(\Oscr,\Oscr)
\]
\begin{theorems}
There is a quasi-isomorphism
\begin{equation}
\label{ref-7.9-81}
\CC_{\mathrm{ab}}(\Pre(\Oscr))\cong \CC_{\mathrm{GS}}(\Oscr)
\end{equation}
\end{theorems}
\begin{proof}
Let the $(P_V)_V$ be as in the proof of the previous
proposition and let $\Oscr!$ be the endomorphism algebra of
the projective generator $P=\coprod_LP_L$ of $\Pre(\Oscr)$. The main
result of \cite{GS2} is the difficult ``Special Cohomology Comparison 
Theorem'':
\[
\RHom_{\Oscr^{\op} \otimes \Oscr}(\Oscr,\Oscr) \cong
\RHom_{\Oscr!^{\op} \otimes \Oscr!}(\Oscr!,\Oscr!)
\]
By Proposition \ref{ref-5.5.1-60} and \eqref{ref-2.5-15} we have
\[
\RHom_{\Oscr!^{\op} \otimes \Oscr!}(\Oscr!,\Oscr!)\cong
\CC_{\mathrm{ab}}(\Pre(\Oscr))
\]
This finishes the proof.
\end{proof}
\subsection{Sheaves of modules over sheaves of rings}
\label{ref-7.3-82} Let $(X,\Oscr)$ be a ringed space.
We prove  \eqref{ref-7.3-75}.
\begin{theorems}\label{ref-7.3.1-83}
Suppose $\bbb$ is a basis of $X$ of \emph{acyclic opens} and put $\BBB =
\UUU({\bbb})$. There is a quasi-isomorphism
\[
\CC(X) \cong \CC_{\mathrm{sh}}({\BBB}).
\]
\end{theorems}
\begin{proof}
Consider the composition
$$j:\BBB \lra \Pre(\ooo) \lra \Mod(\ooo): U \longmapsto i_{U!}\ooo_U$$
where $i_{U!}\ooo_U$ is the  sheafification of $P_U$. Since $\uuu$ is a
basis for the topology, $j$ induces a localization (see for example
\cite{lowen1}). For $U \subset V$, we have
\[
\Ext^i_{\ooo}(i_{U!}\ooo_U,i_{V!}\ooo_V) =
\Ext^i_{\ooo_U}(\ooo_U, (i_{V!}\ooo_V){\mid} U) = H^i(U,\ooo_U)
\]
  and
hence
\begin{enumerate}
\item $\BBB(U,V) \lra \Mod(\ooo)(i_{U!}\ooo_V,i_{V!}\ooo_V)$ is an isomorphism
\item $\Ext^i_{\ooo}(i_{U!}\ooo_V,i_{V!}\ooo_V) = 0$ for $i>0$
\end{enumerate}
So if we endow $\Ob(\BBB)$ with
the censoring relation
\[
(U,V)
\in
\rrr \iff U
\subset V
\]
we see that the result follows from Proposition \ref{ref-5.5.1-60}.
\end{proof}
\begin{remarks}
\label{ref-7.3.2-84} The use of the censoring relation $\Rscr$ is essential
   in the above proof as we have no control over
   $\Ext^i_{\ooo}(i_{U!}\ooo_U,i_{V!}\ooo_V)$ when $U\not\subset V$.
\end{remarks}
\subsection{Constant sheaves}
\label{ref-7.4-85}
In this section we prove \eqref{ref-7.2-74}.
I.e.\ for a topological space $X$ there is an isomorphism
$H\CC^\ast(X,\underline{k})\cong H^\ast(X,\underline{k})$.
The proof is basically a concatenation of some standard facts about
cohomology of presheaves and sheaves.

\medskip

If $\Cscr$ is a small category and $F$ is a presheaf of $k$-modules on
$\Cscr$ then the (presheaf!) cohomology $H^\ast(\Cscr,F)$ of $F$ is
defined as the evaluation at $F$ of the right derived functor
$R^\ast\projlim$ of the inverse limit functor over $\Cscr^{\op}$. It
is well-known that $H^\ast(\Cscr,F)$ can be computed with simplicial
methods \cite{Laudal,Quillen,Roos}. To be more precise put
\[
F^n=\prod_{C_{i_0}
\xrightarrow{f_{i_1}}\cdots \xrightarrow{f_{i_n}} C_{i_n}} F(C_{i_0})
\]
where, as the notation indicates, the product runs over all $n$-tuples
of composable morphisms.
Then $F^\bullet=(F^n)_n$ is a cosimplicial $k$-module and $H^\ast(\Cscr,F)$
is the homology of the associated standard complex
\[
0\r F^0\r F^1\r F^2\cdots
\]
where the differentials are the usual alternating sign linear
combinations of the boundary maps in $F^\bullet$.

Now let  $k\Cscr$ be the $k$-linear
category with the same objects as $\Cscr$ and $\Hom$-sets given by
\[
(k\Cscr)(C,D)=k^{\oplus \Cscr(C,D)}
\]
It is clear that $\Cscr\mapsto k\Cscr$ is the left adjoint to the
forgetful functor from $k$-linear categories to arbitrary
categories. The formula for $F^n$ may be rewritten as
\[
F^n=\prod \Hom_k((k\Cscr)(C_{i_n},C_{i_{n+1}})\otimes_k\cdots \otimes_k
(k\Cscr)(C_{i_0},C_{i_1}),F(C_{i_0}))
\]
where now the  product runs over $n+1$ tuples of objects in $\Cscr$.

To simplify even further assume that $\Cscr$ is a poset.
We make $F$ into a $k\Cscr-k\Cscr$ bimodule in
the
following way
\[
F(C,D)=F(C)
\]
The dependence of $F(C,D)$ on $D$ is as follows: if $f:D\r D'$ is a
map in $\Cscr$ and $\alpha\in k$ then $\alpha f$ is the map from
$F(C,D)=F(C)$ to $F(C,D')=F(C)$
given by multiplication by $\alpha$. With this definition we may
rewrite $F^n$ once again as
\[
F^n=\CC^n(k\Cscr,F)
\]
Now let $X$ be a topological space and put $\Cscr=\Open(X)$, ordered
by inclusion. In addition put $\CCC=k\Cscr$. Let $\underline{k}^p$ be 
the constant presheaf on $X$
with values in $k$ and let $\Oscr=\underline{k}$ be its  sheafification.
We put $\Mod(X)$ for $\Mod(X,\Oscr)=\Mod(X,\underline{k}^p)$,
and similarly $\Pr(X)=\Pr(X,\underline{k}^p)$. Let $\epsilon:\Mod(X)\r
\Pr(X)$ be the inclusion functor.
\begin{convention} To avoid some confusing notations in this section, 
the sections of a (pre)sheaf $\Gscr$
   on an open $U$ will always be denoted by $\Gamma(U,\Gscr)$ and not
   by $\Gscr(U)$.
\end{convention}
If $\Fscr\in \Pr(X)$ then we have
\[
R\projlim\nolimits_{\Cscr^{\op}}(\Fscr)=\projlim\nolimits_{\Cscr^{\op}}(\Fscr)=\Gamma(X,\Fscr)
\]
since $\Cscr^{\op}=\Open(X)^{\op}$ has an initial object $X$. Applying
this to an injective resolution $0\r \Gscr\r I^\cdot$ of an object
$\Gscr$ in $\Mod(X)$
we
find
\[
R\Gamma(X,\Gscr)=\Gamma(X,I^\cdot)=\Gamma(X,\epsilon I^\cdot)=
R\projlim(\epsilon I^\cdot)
\]
and hence
\[
H^\ast(X,\Gscr)=H^\ast(\CC(\CCC,I^\cdot))
\]
where we have suppressed the $\epsilon$.
We will construct an isomorphism
\begin{equation}
\label{ref-7.10-86}
H^\ast(\CC(\CCC,I^\cdot))\cong H^\ast(\CC_{\mathrm{ab}}(\Mod(X)))
\end{equation}
for a specific choice of $I^\cdot$.

Since
$\Mod(\CCC)\cong \Pr(X)$, the functor
\[
j:\CCC\r \Mod(X):U\mapsto i_{U!}\Oscr_U
\]
defines a localization and hence the results from \S\ref{ref-5.1-44}
apply. For $U\in \Ob(\CCC)$ we choose functorial injective resolutions
$U\mapsto E(U)$ of $i_{U!}\Oscr_U$ as in Remark
\ref{ref-5.2.4-51}.

For $U,V\in \Cscr$ put
\[
\bar{\bar{\CCC}}(U,V)=\Hom_{\Mod(X)}(E(U),E(V))
\]
According to Lemma \ref{ref-5.4.2-58} we have
\[
\CC_{\mathrm{ab}}(\Mod(X))=\CC(\CCC,\bar{\bar{\CCC}})
\]
We prove the isomorphism \eqref{ref-7.10-86} for $I^\cdot=E(X)$. To
this end
it is sufficient by Proposition \ref{ref-4.3.1-33} to construct a
quasi-isomorphism between the complexes of $\CCC$-bimodules $E(X)_0$
and $\bar{\bar{\CCC}}_0$. I.e.\
for $U\subset V$ we must construct quasi-isomorphisms between
$E(X)(U,V)$ and $\bar{\bar{\CCC}}(U,V)$ which are natural in
$U,V$. Recall that  in the current setting
\[
E(X)(U,V)=\Gamma(U,E(X))
\]
We have quasi-isomorphisms
\begin{multline*}
\bar{\bar{\CCC}}(U,V)=\Hom_X(E(U),E(V))
\xrightarrow{\cong} \Hom_X(i_{U!}\Oscr_U,E(V))
\xrightarrow{\cong}\\\Gamma(U,E(V))\xrightarrow{\cong}
\Gamma(U,E(X))
\end{multline*}
The last arrow is obtained from the map $E(V)\r E(X)$ which comes from
the map $V\r X$ by functoriality. To see that it is a
quasi-isomorphism note that  $i_{V!}(\Oscr_V)\mid U\cong
\Oscr\mid U$ implies that $E(V)\mid U\r E(X)\mid U$ is a
quasi-isomorphism. Since $E(V)\mid U$ and $E(X)\mid U$ consist of
injectives we obtain indeed a quasi-isomorphism between $\Gamma(U,E(V))$ and
$\Gamma(U,E(X))$.

\subsection{Sheaves of modules over a quasi-compact, separated 
scheme}\label{ref-7.5-87}
In this section we prove \eqref{ref-7.5-77}. Let $X$ be a
quasi-compact separated scheme and let $X = \bigcup_{i=1}^nA_i$ be a
finite affine covering of $X$. For $J \subset I = \{1,2,\dots,n\}$,
put $A_J = \bigcap_{i \in J}A_i$. Each $A_J$ is affine since $X$ is
separated. Put $\aaa = \{A_J\,|\,\varnothing \neq J \subset I\}$.

\begin{theorems}\label{ref-7.5.1-88}
There is a quasi-isomorphism $$\CC(X) \cong
\CC_{\mathrm{sh}}({\UUU(\aaa)}).$$
\end{theorems}
\begin{proof}
   Let $\xxx$ be the collection of all open subsets in $X$ and fix once
   and for all a $k$-cofibrant resolution $ \overline{\UUU(\xxx)}\lra
   \UUU(\xxx)$ with $\overline{\UUU(\xxx)}(U,V) = 0$ if $U$ is not in
   $V$ (recall that we can achieve this by replacing an arbitrary
   $k$-cofibrant resolution $\overline{\UUU(\xxx)}$ by
   $\overline{\UUU(\xxx)}_0$).

All resolutions will be chosen to be
restrictions of
$\overline{\UUU(\xxx)}
\lra
\UUU(\xxx)$.
If $\ccc$, $\ddd$ are collections of opens and $\CCC$ is the $k$-cofibrant
category corresponding to $\ccc$, $\CCC_{\ddd}$ contains the opens $U \in
\ccc$ with $U \subset D$ for some $D \in \ddd$, and $\CCC_D =
\CCC_{\{D\}}$. We write $\ccc \leq \ddd$ (or $\CCC \leq \ddd$) if for every $C
\in
\ccc$, there is a
$D \in \ddd$ with $C \subset D$.

Let $\bbb$ be a basis of affine opens and $\aaa'$ a collection of affine opens
with
$\aaa
\subset
\aaa'
\subset
\bbb$. Let $\AAA$, $\AAA'$ and $\BBB$ denote the corresponding $k$-cofibrant
categories. We will prove that the induced map $\CC(\BBB) \lra
\CC(\AAA'_{\aaa})$ is a quasi-isomorphism.
Taking $\aaa' = \aaa$, the result then follows from Theorem \ref{ref-7.3.1-83}.
The proof goes by induction on the number $n$ of affine
opens in the covering used to produce $\aaa$.

For $n = 1$, the statement follows from the Lemma
\ref{ref-7.5.2-89} below. For arbitrary $n$, put $\aaa_1 = \{A_i\}_{i=2}^n$ and
$\aaa_2 = \{A_1 \cap A_i\}_{i=2}^n$. For any $\CCC \leq \aaa$, we get an exact
sequence of  double complexes (cfr.\ \S\ref{ref-2.4-13})
$$0 \lra \DD(\CCC) \lra \DD(\CCC_{A_1}) \oplus \DD(\CCC_{\aaa_1})
\lra \DD(\CCC_{\aaa_2}) \lra 0$$
since $\CCC(C_{p-1},C_p) \otimes \dots \otimes \CCC(C_0,C_1)$ can be different
from zero only if $C_0 \subset C_1 \subset \dots \subset C_{p-1} \subset C_p$.
Applying this with $\CCC = \BBB, \AAA'$, we obtain a commutative diagram
\[
\xymatrix{0 \ar[r] &\CC(\BBB_{\aaa}) \ar[r]\ar[d]& \CC(\BBB_{A_1}) \oplus
\CC(\BBB_{\aaa_1})
\ar[r] \ar[d]& \CC(\BBB_{\aaa_2}) \ar[r]\ar[d]& 0\\
0 \ar[r] &\CC(\AAA'_{\aaa}) \ar[r]& \CC(\AAA'_{A_1}) \oplus
\CC(\AAA'_{\aaa_1})
\ar[r]& \CC(\AAA'_{\aaa_2}) \ar[r]& 0}
\]
It suffices to prove that the arrow to the left is a quasi-isomorphism. This
follows from the induction hypothesis applied to the three maps to the right.
\end{proof}

\begin{lemmas}\label{ref-7.5.2-89}
Suppose $X$ is affine and let $\ccc$ be a collection of affine subsets with $X
\in \ccc$. Then $\UUU(\{X\}) \lra \UUU(\ccc)$ induces a quasi-isomorphism
$\CC(\ovl{\UUU(\ccc)}) \lra \CC(\ovl{\UUU(\{X\})})$ for restrictions of
$\ovl{\UUU(\xxx)}
\lra \UUU(\xxx)$.
\begin{proof}
Put $\CCC=\ovl{\UUU(\ccc)}$ and $\YYY=\ovl{\UUU(\{X\})}$.  We endow
$\CCC$ with a censoring  relation $(U,V) \in \rrr \iff U \subset V$.
and we use Proposition \ref{ref-4.3.4-38} for the inclusion $\YYY\subset \CCC$.

For $U \subset V \in \ccc$,  we compute
\begin{align*}
   \RHom_{\YYY^{\op}}(\CCC(V,-),\CCC(U,-))
   &\cong \RHom_{\ooo(X)^{\op}}(\ooo(V),\ooo(U))\\
   &\cong \RHom_{\Qch(\ooo_X)}(i_{V,\ast}\ooo_V,i_{U,\ast}\ooo_U)\\
   &\cong \RHom_{\Mod(\ooo_X)}(i_{V,\ast}\ooo_V,i_{U,\ast}\ooo_U)\\
   &\cong \RHom_{\Mod(\ooo_U)}((i_{V,\ast}\ooo_V)|_U, \ooo_U)\\
   &\cong \RHom_{\Mod(\ooo_U)}(\ooo_U,\ooo_U)\\
   &\cong \ooo(U)\\
   &\cong \CCC(U,V)
\end{align*}
where the first line follows from the fact that $\Gamma(X,-)$ defines
an equivalence between $\Mod(\Oscr(X))$ and $\Qch(X)$ and the third
line follows from the separatedness of affine schemes \cite[Appendix
B]{thomasson}.
\end{proof}
\end{lemmas}
\subsection{Computing $\RHom$'s using a covering}
\label{ref-7.6-90}
Let $X$ be a quasi-compact quasi-separated scheme and let $X = \cup_{i
   =1}^nA_i$ be as in the previous section. We use the same associated
notations. Let $\Qch(X)$ be the category of quasi-coherent sheaves on
$X$. In this section we prove that the functor $\Qch(X)\r
\Pre(\Oscr_\Ascr)$ preserves $\RHom$.  The actual reason for this is
that one may show that the simplicial scheme  $S_\bullet$ defined by
 $S_n=\coprod_{i_1\le \cdots \le i_n} A_{i_1}\cap\cdots \cap A_{i_n}$
satisfies ``effective cohomological descent''
\cite[Expos\'e V${}^{\text{bis}}$]{SGA4tome2} for
the obvious map $\epsilon:S_\bullet\r X$, even though it is not quite
a hypercovering in the sense of \cite[Expos\'e V7]{SGA4tome2} (to 
obtain a hypercovering we
need to take all sequences $(i_1, \cdots ,i_n)$ and not just
the ordered ones).  For the convenience of the reader we will give a
direct proof of the preservation of $\RHom$ in our special case.

Before giving the proof let us give a quick sketch. Let
$\epsilon^\ast:\Qch(X)\r \Pr(\Oscr_\Ascr)$ be the exact inclusion
functor. It is easy to see that $\epsilon^\ast$ has a right adjoint
$\epsilon_\ast$, which is some kind of global section functor,
satifying $\epsilon_\ast\epsilon^\ast=\Id$.

We then prove that $\epsilon_\ast$ sends injective objects in
$\Qch(X)$ to acyclic objects for $\epsilon^\ast$. Thus we obtain
$R\epsilon_\ast\circ \epsilon^\ast=\Id$ and hence $\epsilon^\ast$ is
fully faithful for $\RHom$.

\medskip

For convenience, let $\tilde{\Delta}$ be the poset $\{I\,|\, I \subset
\{1,\dots,n\}\}$ ordered by reversed inclusion and let $\Delta$ be its
subposet of all $J \neq \varnothing$. For $I\in \tilde{\Delta}$ put
$A_I=\bigcap_{i\in A} A_i$ (with $A_\varnothing=X$). We have maps
$\Delta \lra \tilde{\Delta} \lra \Open(X):I \longmapsto A_I$ (with
$A_\varnothing=X$) which allow us to consider the restrictions
$\ooo_{\tilde{\Delta}}$ and $\ooo_{\Delta}$ of $\Oscr$. We will think of
$\Pre(\ooo_{\Delta})$ as (equivalent to) the ``category of (presheaf)
objects in the stack of abelian categories $\Qch: \Delta \lra \Cat: I
\longmapsto \Qch(A_I)$''.

In order to abstract the reasoning we will formulate our results in
the following somewhat more general setting.  $\tilde{\sss}$ will be a stack
of Grothendieck categories on $\tilde{\Delta}$ with exact restriction
functors possessing a fully faithful right adjoint, and $\sss$ will be
its restriction to $\Delta$.

For $I \supset J$, we write
$$i^{\ast}_{IJ}: \sss(J) \lra \sss(I)$$
for the exact restriction functor and
$$i_{IJ\ast}: \sss(I) \lra \sss(J).$$
for its fully faithful right adjoint. We will put 
$i^{\ast}_{I\varnothing} = i^{\ast}_{I}$ and
$i_{I\varnothing \ast} = i_{I\ast}$.
  Besides the above mentioned
properties we will also use the following properties
\begin{enumerate}
\item[(C1)] (Base Change) $i^\ast_{I}i_{J\ast}=i_{I\cup
     J,I\ast}i^\ast_{I\cup J,J}$.
\item[(C2)]
If $E
\in \tilde{\sss}(\varnothing)$ is injective, then for every $K \in \Delta$,
$i^{\ast}_KE$ is acyclic for $i_{K\ast}$.
\end{enumerate}
It is not clear to us if (C1) does not follow from the other
properties. (C2) will be only used in the proof of Theorem \ref{ref-7.6.6-96}
which is the main result of this section.
The additional conditions are clearly satisfied for the stacks
$U\mapsto \Mod(U)$ on a finite cover of a ringed space
  and for $U\mapsto \Qch(U)$
on a finite affine cover
of a separated scheme.

A \emph{(presheaf)object in
$\sss$} consists of objects
$(M_I)_{I}$ with
$M_I
\in
\sss(I)$ and maps $(\phi_{IJ})_{I\supset J}$ with $\phi_{IJ}: i^{\ast}_{IJ}M_J
\lra M_I$ for
$I \supset J$ satisfying the obvious compatibilities. Presheaf objects in
$\sss$ and the obvious compatible morphisms constitute a Grothendieck category
$\Pre(\sss)$.

\begin{propositions}\label{ref-7.6.1-91}
The exact functor $$j^{\ast}_I: \Pre(\sss) \lra \sss(K): (M_K)_K \longmapsto
M_I$$ has a right adjoint $j_{I\ast}$ with
\[
(j_{I\ast}M)_J = \begin{cases} i_{IJ\ast}M&\text{if}\,\, I \supset J\\ 0
&\text{otherwise}
\end{cases}
\]
\end{propositions}

\begin{propositions}\label{ref-7.6.2-92}
The exact functor $$x: \tilde{\sss}(\varnothing) \lra \Pre(\sss): M \longmapsto
(i^{\ast}_IM)_I$$
has a right adjoint $$y: \Pre(\sss) \lra \tilde{\sss}(\varnothing): (M_K)_K
\longmapsto \invlim_K(i_{K\ast}M_K).$$
\end{propositions}
\begin{proof}
For $M \in \tilde{\sss}(\varnothing)$, $(N_K)_K \in \Pre(\sss)$, we
have
$$\begin{aligned}
\Hom_{\tilde{\sss}(\varnothing)}(M,y(N_K)_K) &=
\Hom_{\tilde{\sss}(\varnothing)}(M,\projlim_K(i_{K\ast}N_K))\\
&= \projlim_K\Hom_{\tilde{\sss}(\varnothing)}(M,i_{K\ast}N_K)\\
&= \projlim_K\Hom_{\tilde{\sss}(\varnothing)}(i^{\ast}_KM,N_K)\\
&= \Hom_{\Pre(\sss)}(\epsilon^\ast M,(N_K)_K)\qed
\end{aligned}$$
\def\qed{}\end{proof}

For every $K \in \Delta$, we obtain two commutative diagrams:

$$\xymatrix{&{\tilde{\sss}(\varnothing)}\ar[ld]_x \ar[d]^{i^{\ast}_K}\\
{\Pre(\sss)} \ar[r]_{j^{\ast}_K} & \sss(K)}
\hspace{1cm}
\xymatrix{&{\tilde{\sss}(\varnothing)}\\ {\Pre(\sss)} \ar[ru]^y &{\sss(K)}
\ar[u]_{i_{K\ast}} \ar[l]^{j_{K\ast}}}$$
The arrows in the left triangle are exact left adjoints to the corresponding
arrows in the right triangle, which consequently preserve injectives.

Let $\EEE$ denote the full subcategory of $\Pre(\sss)$ spanned by the objects
$j_{K\ast}E$ for $E$ injective in $\sss(K)$ and $K \in \Delta$. Let
$\add(\EEE)$ be spanned by all finite sums of objects in $\EEE$.

\begin{propositions}\label{ref-7.6.3-93}
$\Pre(\sss)$ has enough injectives in $\add(\EEE)$.
\end{propositions}
\begin{proof}
For $M = (M_K)_K$ in $\Pre(\sss)$, we can choose monomorphisms $M_K \lra E_K$
in $\sss(K)$ for each $K$. The corresponding maps $M \lra j_{K\ast}E_K$ yield
a map $$M \lra \bigoplus_Kj_{K\ast}E_K,$$ which is a monomorphism since every
image under $j^{\ast}_K$ is.
\end{proof}
We will now describe
$R\epsilon_\ast$ for certain $N \in \Pre(\sss)$.
To $N = (N_K)_K$ in $\Pre(\sss)$, we associate the following complex $S(N)$
in $C(\tilde{\sss}(\varnothing))$:
$${\prod_pi_{p\ast}N_p} \lra
{\prod_{p<q}i_{pq\ast}N_{pq}} \lra {\dots} \lra {i_{1\dots
n\ast}N_{1\dots n}} \lra 0$$
with ${\prod_pi_{p\ast}N_p}$ in degree zero and with the usual alternating sign
differentials $d^0_N,d^1_N,\dots$.

\begin{propositions}\label{ref-7.6.4-94}
For $N = (N_K)_K$ in $\Pre(\sss)$, we have
\[
\epsilon_\ast N = H^0(S(N))
\]
If for every $K$, $N_K$ is acyclic for $i_{K\ast}$, we have in
$D(\tilde{\sss}(\varnothing))$ $$R\epsilon_\ast N = S(N).$$
\begin{proof}
First note that it is clear from the shape of $\Delta$ that if $(F_K)_{K \in
\Delta}$ is a functor $\Delta \lra \tilde{\sss}(\varnothing)$, then 
$\projlim_{K}F_K =
\projlim_{|K| \in \{1,2\}}F_K$. Hence the first statement follows 
from Proposition
\ref{ref-7.6.2-92}.
For the second statement, first consider $N = j_{K\ast}E \in \EEE$. By
Proposition \ref{ref-7.6.1-91}, the complex $S(N)$ is given by
$${\prod_{p\in K}i_{K\ast}E} \lra
{\prod_{p<q \in K}i_{K\ast}E} \lra {\dots} \lra {i_{K\ast}E} \lra 0.$$
This complex is acyclic, except in degree zero where its homology is
$i_{K\ast}E = \epsilon_\ast N$. It follows that for all $E \in
\add(\EEE)$, $S(E) = \epsilon_\ast E =
R\epsilon_\ast E$.
Now consider $N$ with every $N_K$ acyclic for $i_{K\ast}$.
Take a resolution $N \lra E^{\cdot}$ of $N$ in $\add(\EEE)$. Consider
$S(E^{\cdot})$ as a first quadrant double complex with the complexes
$S(E^i)$ vertical. Looking at columns first, we find a quasi-isomorphism
$\Tot(S(E^{\cdot})) \cong \epsilon_\ast E^{\cdot} \cong R\epsilon_\ast N$.
Since for $E \in \add(\EEE)$, the $E_K$ are obviously acyclic for
$i_{K\ast}$, it follows from our assumption on $N$ that $i_{K\ast}(N_K \lra
E^{\cdot}_K)$ is exact. Hence looking at rows, we find a quasi-isomorphism
$\Tot(S(E^{\cdot})) \cong S(N)$, which finishes the proof.
\end{proof}
\end{propositions}

\begin{propositions}\label{ref-7.6.5-95}
For $M \in \tilde{\sss}(\varnothing)$, we have
\begin{enumerate}
\item $\epsilon_\ast \epsilon^\ast M = M$;
\item $S(\epsilon^\ast M) = M$ in $D(\tilde{\sss}(\varnothing))$;
\end{enumerate}
\end{propositions}
\begin{proof}
To prove the  two statements, it suffices that the complex
$$0 \lra M \lra {\prod_pi_{p\ast}i^{\ast}_pM} \lra
{\prod_{p<q}i_{pq\ast}i^{\ast}_{pq}M} \lra {\dots} \lra {i_{1\dots
n\ast}i^{\ast}_{1\dots n}M} \lra 0$$
is acyclic. This follows from the fact that for every $r \in \{1,\dots n\}$,
the image of the complex under $i^{\ast}_r$ has a contracting chain
homotopy.
\end{proof}
\begin{theorems}\label{ref-7.6.6-96}
The functor $\epsilon^\ast: \tilde{\sss}(\varnothing) \lra 
\Pre(\sss)$ induces a fully
faithful functor $\epsilon^\ast :D^+(\tilde{\sss}(\varnothing))\lra 
D^+(\Pre(\sss))$
\end{theorems}
\begin{proof}
   It is sufficient to prove that $R\epsilon_\ast\circ \epsilon^\ast$
   is the identity on $D^+(\tilde{\sss}(\varnothing))$. To this end it
   is sufficient to prove that $\epsilon_\ast \epsilon^\ast$ is the
   identity and that $\epsilon^\ast$ sends injectives to acyclic objects for
   $\epsilon_\ast$.  The first part is Proposition \ref{ref-7.6.5-95}(1). The
   second part follows from Proposition \ref{ref-7.6.5-95}(1), the extra
   condition (C2) on $\Sscr$ and Proposition \ref{ref-7.6.4-94}.
\end{proof}

\subsection{Quasi-coherent sheaves over a quasi-compact, separated
   scheme}
\label{ref-7.7-97}
We keep the notations of the previous section.  The following theorem
is still true in the general setting exhibited there.
\begin{theorems}\label{ref-7.7.1-98}
There is a quasi-isomorphism 
$$\CC_{\mathrm{ab}}(\tilde{\sss}(\varnothing)) \cong
\CC_{\mathrm{ab}}(\Pre(\sss)).$$
\begin{proof}
   Let $\EEE$ be as before and let $\III(\varnothing)$ be the category
   of injectives in $\tilde{\sss}(\varnothing)$. Consider the
   $\EEE{-}\III(\varnothing)$-bimodule $X$ with
$$X(I,j_{K\ast}E) = \Hom_{\Pre(\sss)}(\epsilon^\ast I,j_{K\ast}E) =
\Hom_{\tilde{\sss}(\varnothing)}(I,i_{K,\ast}E).$$
For $I,J \in \III(\varnothing)$, we compute:
$$\begin{aligned}
\RHom_{\EEE^{\op}}(X(J,-),X(I,-)) &=
\RHom_{\EEE^{\op}}(\Hom_{\Pre(\sss)}(\epsilon^\ast
J,-),\Hom_{\Pre(\sss)}(\epsilon^\ast I,-))\\
&= \RHom_{\Pre(\sss)}(\epsilon^\ast I,\epsilon^\ast J)\\
&= \RHom_{\tilde{\sss}(\varnothing)}(I,J)\\
&= \Hom_{\tilde{\sss}(\varnothing)}(I,J)
\end{aligned}$$
where the second step follows from Lemma \ref{ref-5.3.3-55} and the third step
is Theorem \ref{ref-7.6.6-96}.
For $j_{K\ast}E$ and $j_{L\ast}F$ in $\EEE$, first note that
$$\Hom_{\Pre(\sss)}(j_{K\ast}E, j_{L\ast}F) =
\Hom_{\sss(L)}(j^{\ast}_Lj_{K\ast}E,F),$$ which equals zero unless $L \subset
K$, and equals
\begin{align*}
\Hom_{\sss(L)}(i_{KL\ast}E,F) &=
\Hom_{\tilde{\sss}(\varnothing)}(i_{K\ast}E, i_{L\ast}F)\\
&=\Hom_{\EEE}(i_{K\ast}E, i_{L\ast}F)
\end{align*}
  if
$L
\subset K$ since $i_{L\ast}$ is fully faithful.
For $L \subset K$, we now compute
$$\begin{aligned}
\RHom_{\III(\varnothing)}(\Hom_{\tilde{\sss}(\varnothing)}(-,i_{K\ast}E),
\Hom_{\tilde{\sss}(\varnothing)}(-,i_{L\ast}F))
&= \Hom_{\III(\varnothing)}(i_{K\ast}E,i_{L\ast}F)\\
& = \Hom_{\EEE}(j_{K\ast}E, j_{L\ast}F)\\
&=\Hom_{\CCC}(i_{K\ast}E, i_{L\ast}F)
\end{aligned}$$
We endow $\EEE$ with the
censoring
relation $(j_{K\ast}E, j_{L\ast}F) \in \rrr \iff j_{K\ast}E \neq 0
\neq j_{L\ast}F \,\text{and}\, L \subset K$.
By Proposition \ref{ref-4.3.3-36} and  Corollary \ref{ref-6.7-68} we obtain
the desired quasi-isomorphism.
\end{proof}
\end{theorems}

Specializing to quasi-coherent sheaves on a quasi-compact separated
scheme we get
\begin{corollarys}
\label{ref-7.7.2-99}
There is a quasi-isomorphism $$\CC_{\mathrm{ab}}(\Qch(X)) \cong
\CC_{\mathrm{ab}}(\Pre(\ooo_{\aaa})).$$
\end{corollarys}
\noindent
and in the case of a noetherian separated scheme we get
\begin{corollarys}
\label{ref-7.7.3-100}
There is a quasi-isomorphism $$\CC_{\mathrm{ab}}(\coh(X)) \cong
\CC_{\mathrm{ab}}(\Pre(\ooo_{\aaa})).$$
\end{corollarys}
\noindent
using Corollary \ref{ref-6.8-69}, since $\Qch(X)=\Ind\coh(X)$.

\begin{remarks}
By combining the above results we obtain an isomorphism in $\Ho(B_\infty)$
\begin{equation}
\label{ref-7.11-101}
\CC_{\text{ab}}(\Qch(X))\cong \CC_{\text{ab}}(\Mod(X))
\end{equation}
  for a
quasi-compact separated scheme $X$, but our proof
of this fact is far from straightforward and goes through the
auxilliary categories $\Pre(\Oscr_\Ascr)$ and $\Pre(\Oscr_\Bscr)$. It
would be interesting to see if a more direct proof could be obtained. 
\end{remarks}
\subsection{Relation to Swan's definition}

In \cite{swan}, Richard G. Swan defined the Hochschild cohomology of a
separated scheme $X$ to be
$$\Ext^i_{X\times X}(\ooo_D,\ooo_D)$$
where $\ooo_D = \delta_{\ast}\ooo_X$ for the diagonal map $\delta: X \lra X
\times X$.
Put $\CC_{\mathrm{Swan}}(X) = \RHom_{X \times X}(\ooo_D,\ooo_D)$.

We prove that Swan's definition coincides with ours. As was already
mentioned this in the finite type case could be deduced from \cite[\S
3]{swan} and the above results.
\begin{theorems}
\label{ref-7.1-102}
   Let $X$ be a quasi-compact, separated scheme over a field $k$ with
   an affine covering $\Ascr$ given by $X = \cup_{i =1}^nA_i$. There is a
   quasi-isomorphism
$$\CC_{\mathrm{Swan}}(X) \cong \CC_{\mathrm{GS}}(\ooo_{\aaa}).$$
\begin{proof}
Consider the factorization of
$\delta: X \lra X \times X$
over $\delta': X \lra X'$ where $X' = \cup_{i =1}^n
A_i \times A_i$ is the open subscheme of $X
\times X$ with $\ooo_{X'}(U) = \ooo_{X \times X}(U)$ and in particular
$\ooo_{X'}(A_i \times A_i) = \ooo_X(A_i) \otimes \ooo_X(A_i)$.
Put $\ooo_{D'} = \delta'_{\ast}\ooo_X$. In particular, $\ooo_{D'}(A_i \times
A_i) = \ooo_X(A_i)$. If we identify the collection $\aaa'$ associated to the
covering $\cup_{i =1}^n
A_i \times A_i$ of $X'$ with $\aaa$, we have ${\ooo_{X'}}_{\aaa'} =
{\ooo_X}_{\aaa} \otimes {\ooo_X}_{\aaa}$ and ${\ooo_{D'}}_{\aaa'} =
{\ooo_X}_{\aaa}$. Since $X$ is separated, we have $\mathrm{sup}(\ooo_D)
\subset \delta(X) \subset X'$, hence we may compute
\begin{align*}
\RHom_{X \times X}(\ooo_D,\ooo_D) &= \RHom_{X'}(\ooo_{D'},\ooo_{D'})\\
&= \RHom_{\Qch(X')}(\ooo_{D'},\ooo_{D'})\\
&= \RHom_{\Pre(\ooo_{\aaa} \otimes
\ooo_{\aaa})}({\ooo_{\aaa}},{\ooo_{\aaa}})
\end{align*}
where we have used that $X$ is quasi-compact, separated in the second step and
we have used Theorem
\ref{ref-7.6.6-96} in the last step.
\end{proof}
\end{theorems}
\begin{corollarys} We have
\[
\CC(X)\cong\CC_{\mathrm{Swan}}(X)
\]
\end{corollarys}
\begin{proof}
We have
\[
\CC_{\mathrm{Swan}}(X)\cong \CC_{GS}(\Oscr_\Ascr)\cong 
\CC(\Pre(\Oscr_\AAA))\cong
\CC(\frak{u}(\AAA)) \cong\CC(X)
\]
The first isomorphism is Theorem \ref{ref-7.1-102}, the second isomorphism
is \eqref{ref-7.9-81}, the third isomorphism is \eqref{ref-7.8-80}
and the fourth isomorphism is Theorem \ref{ref-7.5.1-88}.
\end{proof}
\subsection{The Mayer-Vietoris sequence}
\label{ref-7.9-103}
Let $(X,\Oscr)$ be a ringed space and let $X=U\cup V$. In this section
we prove \eqref{ref-7.1-73}. If $W_1\subset W_2$ are open embeddings of
  ringed spaces  then the pushforward functor
$i_{W_2,W_1\ast}$ is fully faithful and preserves injectives. Hence it
induces a restriction map  $\CC(W_2)\r
\CC(W_1)$
(see Remark \ref{ref-4.2.3-31}). It is clear that the restriction map is
compatible with compositions. So $\CC(-)$ defines a contravariant
functor on open embeddings of ringed spaces.
\begin{theorems}
\label{ref-7.9.1-104}
There is an exact triangle of
complexes
\begin{equation}
\label{ref-7.12-105}
\CC(X)\r \CC(U)\oplus \CC(V)\r \CC(U\cap V)\r
\end{equation}
where the
maps are the restriction maps defined above (in
particular they are $B_\infty$-maps).
\end{theorems}
Taking homology in \eqref{ref-7.12-105} yields the Mayer-Vietoris sequence.
\begin{proof}
We use notations as in \S\ref{ref-7.6-90}. Put $U_1=U$, $U_2=V$ and
$\Delta=\{\{1\},\{2\},\{1,2\}\}$. Let
$\Sscr$ be the stack of abelian categories $\Mod(-)$ associated to the
covering $X=U_1\cup U_2$. Let $\EEE$ be the full subcategory of
$\Pre(\Sscr)$ consisting of $j_{K\ast}E$ for $E$ injective in $\sss(K)$
and $K \in \Delta$.

Let $\EEE_I=\Inj\Mod(U_I)$ for $I\subset \{1,2\}$. As mentioned
above the
functors $i_{I,\varnothing\ast}:\EEE_I\r \EEE_\varnothing$ are fully
faithful.
Choose a $k$-cofibrant resolution $\bar{\EEE}_\varnothing\r
\EEE_\varnothing$ and let $\bar{\EEE}_I\r \EEE_I$ be the restrictions of
this resolution. Let $\bar{\EEE}$ be the category with $\Ob(\bar{\EEE}) =
\coprod_{J \in \Delta}\Ob(\EEE_{J})$ and
$\Hom$-sets  between $E\in \EEE_I$ and $F\in \EEE_J$ given by
\[
\bar{\EEE}(E,F)=
\begin{cases}
\bar{\EEE}_\varnothing(E,F)&\text{if $J\subset I$}\\
0&\text{otherwise}
\end{cases}
\]
Then $\bar{\EEE}$ is a $k$-cofibrant resolution of $\EEE$.

Let $\bar{\EEE}^I$ be the full-subcategory of $\bar{\EEE}$ consisting
of objects in $\bar{\EEE}_J$ with $I\subset J$. We have
$\bar{\EEE}^{\{1,2\}}=\bar{\EEE}_{\{1,2\}}$ and $\bar{\EEE}^{\{j\}}$ is the
     arrow category $(\EEE_{\{1,2\}}\r \EEE_{\{j\}})$ where the arrow
     is the inclusion.

Then by using
restriction maps we obtain a commutative diagram
of complexes
\[
\begin{CD}
0@>>> \CC(\bar{\EEE})@>>> \CC(\bar{\EEE}^{\{1\}})\oplus \CC(\bar{\EEE}^{\{2\}})
@>>> \CC(\bar{\EEE}^{\{1,2\}})@>>> 0\\
@. @. @VVV @VVV \\
@. @. \CC(\bar{\EEE}_{\{1\}})\oplus \CC(\bar{\EEE}_{\{2\}})
@>>> \CC(\bar{\EEE}_{\{1,2\}})@>>> 0
\end{CD}
\]
The vertical maps are isomorphisms by Remark \ref{ref-4.1.4-29} and the top
row is an exact sequence of complexes. Thus we obtain a triangle
\[
\CC(\bar{\EEE})\r \CC(\bar{\EEE}_{\{1\}})\oplus \CC(\bar{\EEE}_{\{2\}})
\r \CC(\bar{\EEE}_{\{1,2\}})\r
\]
  By Theorem
\ref{ref-7.7.1-98}, Proposition \ref{ref-7.6.3-93} and Corollary
\ref{ref-6.7-68} we have $\CC(X)\cong\CC(\bar{\EEE}_\varnothing)\cong
\CC(\bar{\EEE})$.  This means that we are almost done, except for the
fact that we still need to show that the composition
\[
\CC(\bar{\EEE}_\varnothing)\xrightarrow{\cong}\CC(\bar{\EEE})\r
\CC(\bar{\EEE}_{\{j\}})
\]
for $j=1,2$ is the restriction map. To this end we recall that the
isomorphism
$\CC(\bar{\EEE}_\varnothing)\xrightarrow{\cong}\CC(\bar{\EEE})$ was
constructed in the proof of Theorem \ref{ref-7.7.1-98} using the
$\bar{\EEE}-\bar{\EEE}_\varnothing$ bimodule $X$ where
\[
X(E,F)=\Hom_{\Pre(\Sscr)}(\epsilon^\ast E,F)\qquad \text{for $E\in
   \bar{\EEE}_\varnothing$, $F\in \bar{\EEE}$}
\]
We have a commutative diagram of inclusions:
\[
\begin{CD}
   \bar{\EEE}_\varnothing @>\cong>>
   (\bar{\EEE}_\varnothing\xrightarrow{X}
   \bar{\EEE}) @<\cong<< \bar{\EEE}\\
   @| @AAA @AAA \\
\bar{\EEE}_\varnothing @>>i_{\bar{\EEE}_\varnothing}>
(\bar{\EEE}_\varnothing\xrightarrow{X_j}
   \bar{\EEE}_{\{j\}}) @<<< \bar{\EEE}_{\{j\}}
\end{CD}
\]
where ``$\cong$'' means inducing an isomorphism in Hochschild
cohomology. Here $X_j$ is the restriction of $X$ to an
$\bar{\EEE}_{\{j\}}-\bar{\EEE}_\varnothing$ bimodule. It is easy to see that
$X_i$ is the bimodule associated to the inclusion
$i_{\{j\},\varnothing \ast}:\bar{\EEE}_{\{j\}}\lra \bar{\EEE}_{\varnothing}$
(as introduced in Theorem
\ref{ref-4.1.2-26}.2). Thus in the above diagram
$i_{\bar{\EEE}_{\varnothing}}$ induces an isomorphism on Hochschild cohomology
and we obtain a commutative diagram:
\[
\begin{CD}
  \CC( \bar{\EEE}_\varnothing) @>\cong>>
\CC( \bar{\EEE})\\
   @| @VVV \\
\CC(\bar{\EEE}_\varnothing) @>>>  \CC(\bar{\EEE}_{\{j\}})
\end{CD}
\]
where now the lower and the rightmost maps are restriction maps. This
finishes the proof.
\end{proof}


\def\cprime{$'$}
\providecommand{\bysame}{\leavevmode\hbox to3em{\hrulefill}\thinspace}

\end{document}